\documentclass[11pt,a4paper,reqno]{amsart}
  
\usepackage{amsfonts}
\usepackage{amsmath,amssymb,amsthm,amsxtra}
\usepackage{float}
\usepackage{mathtools}
\usepackage[colorlinks, linkcolor=blue!50,anchorcolor=Periwinkle,
    citecolor=blue!72,urlcolor=cyan, bookmarksopen,bookmarksdepth=2]{hyperref}
\usepackage[usenames,dvipsnames]{xcolor}
\usepackage{enumitem}
\usepackage{geometry,array} 
\usepackage{array} 
\usepackage{graphicx}
\usepackage{subfigure}
\usepackage{bookmark}
\usepackage{todonotes}
\usepackage{tikz}
\usepackage{url}
\usetikzlibrary{matrix,positioning,decorations.markings,arrows,decorations.pathmorphing,backgrounds,fit,positioning,shapes.symbols,chains,shadings,fadings,calc}
\tikzset{->-/.style={decoration={  markings,  mark=at position #1 with
    {\arrow{>}}},postaction={decorate}}}
\tikzset{-<-/.style={decoration={  markings,  mark=at position #1 with
    {\arrow{<}}},postaction={decorate}}}

\usepackage[all]{xy}
\def\XX{\mathbb{X}}

\def\QStab{\operatorname{QStab}}

\theoremstyle{plain}
\newtheorem{theorem}{Theorem}[section]

\newtheorem{lemma}[theorem]{Lemma}
\newtheorem{corollary}[theorem]{Corollary}
\newtheorem{proposition}[theorem]{Proposition}

\theoremstyle{definition}
\newtheorem{definition}[theorem]{Definition}

\newtheorem{remark}[theorem]{Remark}

\numberwithin{equation}{section}

\newtheorem{construction}[theorem]{Construction}

\def\hh{\mathcal}

\def\<{\langle}
\def\>{\rangle}
\def\ZZ{\mathbb{Z}}

\def\R{\mathbb{R}}
\def\RR{\R}
\usepackage{verbatim}
\def\cA{\mathcal{A}}

\def\cD{\mathcal{D}}
\def\cE{\mathcal{E}}

\def\cO{\mathcal{O}}

\def\vvt{\tilde{\mathrm{v}}} % numerical Chern character
\def\vv{{\mathrm{v}}} % reduced Chern character

\newcommand{\ch}{\mathrm{ch}}
\newcommand{\glt}{\widetilde{\mathrm{GL}^+(2,\mathbb R)}}

\newcommand{\cccp}{\{1,\frac{\mathrm{ch}_1}{\mathrm{ch}_0},\frac{\mathrm{ch}_2}{\mathrm{ch}_0}\}\text{-plane}}

\makeatletter
\newcommand{\xRrightarrow}[2][]{\ext@arrow 0359\Rrightarrowfill@{#1}{#2}}
\newcommand{\Rrightarrowfill@}{\arrowfill@\equiv\equiv\Rrightarrow}
\newcommand{\xLleftarrow}[2][]{\ext@arrow 3095\Lleftarrowfill@{#1}{#2}}
\newcommand{\Lleftarrowfill@}{\arrowfill@\Lleftarrow\equiv\equiv}
\makeatother

\def\Aut{\operatorname{Aut}}

\def\Hom{\operatorname{Hom}}

\def\Stab{\operatorname{Stab}}

\def\Stap{\operatorname{Stab}^\dag}

\def\deg{\operatorname{deg}}
\renewcommand{\k}{\mathbf{k}}
\renewcommand{\mod}{\operatorname{mod}}

\newcommand{\Cone}{\operatorname{Cone}}

\newcommand{\D}{\operatorname{\hh{D}}}

\def\coh{\operatorname{Coh}}

\def\PP{\mathbb{P}}

\def\gldim{\operatorname{gldim}}
\def\gd{\operatorname{gd}}

\def\iso{\xrightarrow{\sim}}
\def\sli{\mathcal{P}}

\def\CC{\mathbb{C}}

\def\O{\mathcal{O}}

\newcommand{\sslash}{\mathbin{/\mkern-6mu/}}

\newcommand\Geo{\operatorname{Geo}}
\newcommand\Alg{\operatorname{Alg}}

\newcommand{\OO}{\mathcal{O}}
\def\ch{\mathrm{ch}}

\title[Contractibility of $\Stab(\PP^2)$ via $\gldim$]
{Contractibility of space of stability conditions on the projective plane via global dimension function}
\author{Yu-Wei Fan}
\address{Y-W. F:
Yau Mathematical Sciences Center\\
Tsinghua University\\
Beijing, China}
\email{yuweifanx@gmail.com}
\author{Chunyi Li}
\address{C. L:
Mathematics Institute\\
University of Warwick\\
UK}
\email{c.li.25@warwick.ac.uk}

\author{Wanmin Liu}
\address{W. L:
Department of Mathematics\\
Uppsala University\\
Sweden}
\email{wanminliu@gmail.com}
\urladdr{https://wanminliu.github.io/}

\author{Yu Qiu}
\address{Y. Q:
Yau Mathematical Sciences Center and Department of Mathematical Sciences\\
Tsinghua University\\
100084 Beijing, China.
    \& 
Beijing Institute of Mathematical Sciences and Applications\\
Yanqi Lake, Beijing, China}
\email{yu.qiu@bath.edu}

\begin{document}
%=========================================================
%=========================================================
\begin{abstract}
We compute the global dimension function $\mathrm{gldim}$ on the principal component $\mathrm{Stab}^{\dag}(\mathbb{P}^2)$ of the space of Bridgeland stability conditions on $\mathbb{P}^2$. It admits $2$ as the minimum value and the preimage $\mathrm{gldim}^{-1}(2)$ is contained in the closure $\overline{\mathrm{Stab}^{\mathrm{Geo}}(\mathbb{P}^2)}$ of the subspace consisting of geometric stability conditions. We show that $\mathrm{gldim}^{-1}[2,x)$ contracts to $\mathrm{gldim}^{-1}(2)$ for any real number $x\geq 2$ and that $\mathrm{gldim}^{-1}(2)$ is contractible.

    \vskip .3cm
    {\parindent =0pt
    \it Key words.}
    Bridgeland stability conditions, Coherent sheaves, Contractibility, Global dimension function.
    
    \vskip .3cm
    {\parindent =0pt
    \it 2020 Mathematics Subject Classification.}
    14F08 (18G80, 32Q55).

\end{abstract}
\maketitle

%\setcounter{tocdepth}{1} % hide subsections
%\tableofcontents % Table of Contents

%=========================================================
\section{Introduction}
%=========================================================
\subsection{Stability conditions}
%=========================================================
The notion of stability conditions on triangulated categories was introduced by Bridgeland \cite{Bri07},
with motivation coming from string theory and mirror symmetry.
Let $\D$ be a triangulated category and $K_{\mathrm{num}}(\D)$ be its numerical Grothendieck group. A stability condition $\sigma=(Z,\hh{P})$ consists of a central charge
$Z\in\Hom(K(\D),\CC)$ and a slicing $\hh{P}$, which is an $\RR$-collection of t-structures on $\D$. In this paper, we denote $\Stab(\D)$ as
the stability manifold of stability conditions with support property with respect to $K_{\mathrm{num}}(\D)$. By the seminal result in \cite{Bri07}, when $K_{\mathrm{num}}(\D)$ is of finite rank, the space $\Stab(\D)$
is a complex manifold with local coordinate given by the central charge.  The original conjecture \cite[Conjecture 1.2]{Bri08} in the K3 surface case is that $\Stab(\D)$ has a connected component $\Stab^\dag(\D)$ which is simply-connected and preserved by the autoequivalence group of $\D$.
A more ambitious conjecture expects that the stability manifold $\Stab(\D)$ is  contractible in general.
The contractibility is confirmed in a couple of examples at least for the principal component of the space, namely:
\begin{itemize}
  \item The smooth curves case in \cite{Oka06,Bri07,Mac07}.
  \item The K3 surfaces with Picard rank one in \cite{Bri08,BB17}.
  \item The local $\PP^1$ in \cite{IUU10}; the local $\PP^2$ in \cite{BM11}.
  \item The projective plane $\PP^2$  in \cite{Li:spaceP2}.
  \item The Abelian surfaces  in \cite{Bri08} and Abelian threefolds with Picard rank one in \cite{BMS16}.
  \item The finite type (connected) component $\Stab_0$ in \cite{QW18},
  where the heart of any stability conditions in $\Stab_0$
  is a length category with finite many torsion pairs.
  The key examples are (Calabi--Yau) ADE Dynkin quiver case
  and new classes of examples are studied in \cite{AW22}.
  \item The Calabi--Yau-$3$ affine type $A$ case in \cite{Qiu16}.
  \item The acyclic triangular quiver case in \cite{DK16}.
  \item The wild Kronecker quiver case in \cite{DK19}.
\end{itemize}
The proofs in each case are quite different.

%=========================================================
\subsection{Global dimension functions}
%=========================================================
Recently, Ikeda and the fourth-named author \cite{Qiu:gldim,IQ18a} introduce the global dimension  function $\gldim$ on $\Stab(\D)$, namely:
\begin{equation}\label{def:gldim domain}
\gldim\colon \Stab(\D)\to\RR_{\ge0}\cup\{+\infty\},
\end{equation}
which is given by
\begin{equation}\label{def:gldim value}
\gldim \sigma = \gldim\hh{P} \coloneqq \sup\{ \phi_2-\phi_1 \mid
    \Hom(\hh{P}(\phi_1),\hh{P}(\phi_2))\neq0\}.
\end{equation}
Such a function is continuous and invariant under the natural left action by $\Aut(\D)$ and the right action of $\CC$, and thus descends to a continuous function
\begin{equation}\label{def:gldim decent}
\gldim\colon \Aut(\D)\backslash\Stab(\D)/\CC\to\RR_{\ge0}\cup\{+\infty\}.
\end{equation}
The philosophy in \cite{Qiu:gldim} is as follows:
\begin{itemize}
  \item[(i)] The infimum of $\gldim$ on $\Stab(\D)$ (or the principal component of it)
  should be considered as the global dimension $\gd\D$ of the category $\D$.
  \item[(ii)] If the subspace $\gldim^{-1}(\gd\D)$ is non-empty,
  then it is contractible.
  Moreover, the preimage $\gldim^{-1}([\gd\D,x))$ contracts to $\gldim^{-1}(\gd\D)$
  for any real number $\gd\D<x$.
  \item[(iii)] If $\gldim^{-1}(\gd\D)$ is empty,
  then the preimage $\gldim^{-1}(\gd\D,x)$ contracts to $\gldim^{-1}(\gd\D,y)$
  for any real number $\gd\D<y<x$.
\end{itemize}
Note that for a Calabi--Yau category, the global dimension function is constant. If the global dimension function $\gldim$ is not constant, it sheds some lights on why $\Stab(\D)$ should be contractible.

The theme in \cite{IQ18a} is to $q$-deform stability conditions.
More precisely, given a Calabi--Yau-$\infty$ category $\D_\infty$ (e.g. bounded derived category of $\PP^2$),
the corresponding Calabi--Yau-$N$ category $\D_N$ (e.g. local $\PP^2$ for $\PP^2$ and $N=3$) can be obtained by Calabi--Yau-$\XX$ completing $\D_\infty$ to $\D_\XX$
and specializing $\XX$ to be $N$, in other words, taking the orbit category $\D_N=\D_\XX\sslash[\XX-N]$. Under this procedure,  a stability condition $\sigma$ on  $\D_\infty$ such that
\begin{gather}\label{eq:<}
    \gldim\sigma\le N-1
\end{gather}
induces a stability condition on $\D_N$ via  $q$-stability conditions on $\D_\XX$.
We will discuss such inducing in Section~\ref{sec:inducing} for the example from $\PP^2$ to local $\PP^2$ where $N=3$.

%=========================================================
\subsection{The projective plane case}
%=========================================================
In this paper, we study the case of the projective plane $\PP^2$ for the above conjectures/philosophy. The main result is a computation of the global dimension function for the principal component $\Stab^\dag(\PP^2)$ (i.e. the connected component which contains geometric stability conditions, where a stability condition $\sigma\in \Stab(\PP^2)$ is called geometric if all skyscraper sheaves are $\sigma$-stable of the same phase). Details are in Propositions \ref{prop:gldiminParabola} and \ref{prop:glonthestab}. Based on the computation of $\gldim$, we prove the following theorem.

\begin{theorem}[Corollary~\ref{cor:contra2} and Theorem~\ref{thm:3}]\label{thm:1}
Consider the function
$$
    \gldim\colon\Stap(\PP^2)\to\RR_{\geq 0}
$$
on the principal component $\Stap(\PP^2)$
of the space of stability conditions on the bounded derived category $\D=\D^b(\coh \PP^2)$
of coherent sheaves on $\PP^2$.
Then
\begin{itemize}
\item $\gd \D=2$ and $\gldim\Stap(\PP^2)=[2,\infty)$,
\item the subspace $\gldim^{-1}[2,x)$ contracts to $\gldim^{-1}(2)$, for any $x\ge2$,
\item the subspace $\gldim^{-1}(2)$ is contractible and is contained in $\overline{\Stab^{\Geo}(\PP^2)}$,
where $\Stab^{\Geo}(\PP^2)$ consists of geometric stability conditions.
\end{itemize}
\end{theorem}

The contractibility of $\Stap(\PP^2)$ is already proved by the second-named author \cite{Li:spaceP2}. The new approach here shows how this stability manifold contracts along the values of the global dimension function.

%=========================================================
\subsection{Topological Fukaya case}
%=========================================================
In the parallel work \cite{Qiu:flow}, we use the same philosophy to study the contractibility of the space of stability conditions on the topological Fukaya category of a graded marked surface. 
We prove a slightly weaker version of the corresponding 
Theorem~\ref{thm:1}, that $\gldim$ induces the contractible flow except for certain possible critical values.

We hope that these works will shed lights on how this philosophy would apply to other cases.

\subsection*{Acknowledgements}
%=========================================================
C. Li is supported by  the Royal Society URF$\backslash$R1$\backslash$201129 ``Stability condition and application in algebraic geometry" and the Leverhulme Trust ECF-2017-222.
W. Liu is supported by a grant from the Knut and Alice Wallenberg Foundation. He would like to thank Tobias Ekholm and Ludmil Katzarkov for comments.
Y. Qiu is supported by National Key R\&D Program of China (No. 2020YFA0713000), Beijing Natural Science Foundation (Grant No. Z180003) and National Natural Science Foundation of China (Grant No. 12031007).

%=========================================================
\section{Preliminaries}
%=========================================================

%=========================================================
\subsection{The category}\label{sec:func}
%=========================================================
In this paper, we let $\PP^2$ be the projective plane over the complex number field. We write
\begin{gather}\label{eq:D P2}
    \D_\infty(\PP^2)\coloneqq \D^b(\PP^2)=\D^b(\coh\PP^2)
\end{gather}
for the bounded derived category of coherent sheaves on $\PP^2$. Due to the well-known result by Be{\u\i}linson \cite{Bei83}, we have the equivalent description $\D_\infty(\PP^2)\cong\D^b(\k Q/R)$,
where $(Q,R)$ is the quiver
\[
\begin{tikzpicture}[xscale=0.5,
  arrow/.style={<-,>=stealth},
  equalto/.style={double,double distance=2pt},
  mapto/.style={|->}]
\node (x2) at (6,-1){};
\node (x1) at (2,-1){};
\node (x3) at (-2,-1){};
\draw[font=\scriptsize](0,-1)node[above]{$x_1,y_1,z_1$}
    (4,-1)node[above]{$x_2,y_2,z_2$};
  \node at (x1){$2$};
  \node at (x2){$3$};
  \node at (x3){$1$};
\foreach \n/\m in {1/3,2/1}
    {\draw[arrow] (x\n.150) to (x\m.30);
     \draw[arrow] (x\n.-150) to (x\m.-30);
     \draw[arrow] (x\n) to (x\m);}
\end{tikzpicture}
\]
with commutative relations
$$  a_1 b_2=b_1 a_2,\quad a,b\in\{x,y,z\}.$$

\def\Serre{\mathbb{S}} 
The Serre functor on $\D_\infty(\PP^2)$ is given by (see \cite{BK90} or \cite{Huy06})
\[
    \Serre=\Serre_{\PP^2} \coloneqq (-)\otimes \omega_{{\PP}^2}[2]=(-)\otimes \mathcal{O}_{{\PP}^2}(-3)[2].
\]

An object $E\in \D_\infty(\PP^2)$ is called exceptional if $\Hom(E,E[i])=0$ for $i\neq 0$ and $\Hom(E,E)=\CC$. The right and left mutations of an object $F$ with respect to an exceptional object $E$ are defined by
\begin{align}
    \mathsf{R}_{E}(F) &\coloneqq \Cone\left(F\xrightarrow{\mathrm{ev}}E\otimes \Hom(F,E)^*\right)[-1], \label{eq:rightmutation}\\
    \mathsf{L}_{E}(F) &\coloneqq \Cone\left(E\otimes \Hom(E,F) \xrightarrow{\mathrm{ev}}F\right). \label{eq:leftmutation}
\end{align}

%=========================================================
\subsection{An affine plane}
%=========================================================
Let $\D=\D_\infty(\PP^2)$. Let $H$ be the hyperplane divisor of $\PP^2$. For $E\in \D$, we identity the Chern character $\ch(E)$ with the triple of numbers $$\vvt(E)=(\ch_0(E),\ch_1(E).H,\ch_2(E)).$$
When we say the \emph{point} $E$ (or the \emph{point} $\vvt(E)$), we mean the point in the real projective space $\PP(\mathbb{R}^3)$ with homogeneous coordinate $[\ch_0(E),\ch_1(E).H,\ch_2(E)]$. We call the locus $\ch_0=0$ as the line at infinity and its complement as the affine $\cccp$. Moreover, we always
assume that the $\frac{\ch_1}{\ch_0}$-axis is horizontal and the  $\frac{\ch_2}{\ch_0}$-axis is vertical. If $\ch_0(E)\neq 0$, the \emph{reduced character} of $E$ corresponds to the point
\begin{equation}
\vv(E)\coloneqq(1,s(E),q(E)), \quad \text{ with } s(E)\coloneqq\frac{\ch_1(E).H}{\ch_0(E)},\quad q(E)\coloneqq\frac{\ch_2(E)}{\ch_0(E)},
\label{eq:s-q-coordinate}
\end{equation}
in the affine $\cccp$. In particular $\vv(E)=\vv(E[n])$, i.e. $E$ and its any shift $E[n]$ will be the same point in the  $\cccp$.

The $\cccp$ provides a playground for studying both geometric and algebraic stability conditions in the following part of the paper.

%=========================================================
\subsection{Stability conditions}
%=========================================================
A {\it stability condition} $\sigma = (Z, \sli)$ on $\D$ consists of
a group homomorphism $Z \colon K(\D) \to \CC$ called the {\it central charge} and
a family of full additive subcategories $\sli (\phi) \subset \D$ for $\phi \in \R$
called the {\it slicing} satisfying certain conditions.
We refer to \cite{Bri07} and the lecture notes \cite[Definition 5.8]{MS17} for the details.
Nonzero objects in $\sli(\phi)$ are called {\it semistable of phase $\phi$} and simple objects
in $\sli(\phi)$ are called {\it stable of phase $\phi$}.
For semistable object $E\in\sli(\phi)$,
denote by $\phi_\sigma(E)=\phi$ its \emph{phase}.

Let $\D=\D_\infty(\PP^2)$ and
$$\Stab(\PP^2)\coloneqq\Stab(\D_\infty(\PP^2))$$
be the space of stability conditions on $\D_\infty(\PP^2)$. 

A stability condition $\sigma\in \Stab(\PP^2)$ is called \emph{geometric} if all skyscraper sheaves are $\sigma$-stable of the same phase. We denote the set of all geometric stability conditions by $\Stab^{\Geo}(\PP^2)$. 

Let us briefly recall the construction of geometric stability conditions. There is a fractal curve $\mathrm{C}_{\mathrm{LP}}$, the so called Le Potier curve, and a region $\Geo_{\mathrm{LP}}$ in the  $\cccp$, as in Definition~\ref{def:LPcurve}. For each $(1,s,q)\in \Geo_{\mathrm{LP}}$, one can associate a geometric stability condition $\sigma_{s,q}=(Z_{s,q}, \hh{P}_{s,q})$ as follows. The central charge $Z_{s,q}$ is given by
\begin{equation} \label{eq:central charge}
Z_{s,q}(E) \coloneqq (-\ch_2(E)+q\cdot \ch_0(E))+i (\ch_1(E).H-s\cdot \ch_0(E)), \quad \text{ for } E\in \D.
\end{equation}
Denote $H$-slope of coherent sheaves by $\frac{\ch_1(-).H}{\ch_0(-)}$. We make a convention that $H$-slope of a torsion sheaf is $+\infty$. The heart $\hh{P}_{s,q}((0,1])$ is the tilting
$$
    \coh_{\# s}\coloneqq\langle \coh_{\leq s}[1], \coh_{>s}\rangle,
$$
where $\coh_{\leq s}$ (resp. $\coh_{>s}$) is the subcategory of $\coh{(\PP^2)}$ generated by $H$-slope semistable sheaves of slope $\leq s$ (resp. $>s$) by extension.
The slicing for $\phi\in (0,1]$ is defined by
$$
\hh{P}_{s,q}(\phi)=\{E\in \coh_{\# s} \, \, |\, \, E \text{ is } \sigma_{s,q}\text{ -semistable of phase } \phi \}\cup \{0\}.
$$
For general $\phi\in\mathbb{R}$, we have $\hh{P}_{s,q}(\phi+1)=\hh{P}_{s,q}(\phi)[1]$.

The $\glt$ acts freely on $ \Stab^{\Geo}(\PP^2)$ (\cite[Definition 1.4, Corollary 1.15]{Li:spaceP2}) with quotient
$$
\Stab^{\Geo}(\PP^2) / \glt \cong \Geo_{\mathrm{LP}}.
$$

We refer to Section~\ref{sec:alg} for the definition of algebraic stability conditions $\Stab^{\Alg}(\PP^2)$.
Let $\Stap(\PP^2)$ be the connected component in $\Stab(\PP^2)$ which contains the geometric stability conditions. It is still a \emph{conjecture} that $\Stab(\PP^2)=\Stap(\PP^2)$. The second-named author \cite{Li:spaceP2} shows that
\begin{equation}
   \Stap(\PP^2)=\Stab^{\Geo}(\PP^2)\bigcup\Stab^{\Alg}(\PP^2)
   \label{eq:stabmfd}
\end{equation}
and it is contractible. In the following sections, we will compute the global dimension function $\gldim$ on $\Stap(\PP^2)$ and show that the contraction is along the value of $\gldim$.

%=========================================================
\section{Geometric stability conditions in the parabolic region}
%=========================================================

Let $\D=\D_\infty(\PP^2)$. For $a\in \mathbb{R}$, denote by $\Delta_a$ the parabola in the  $\cccp$:
$$
\Delta_a\coloneqq\Big\{(1,s,q)\in \cccp\; | \; \frac{1}{2}s^2-q=a\Big\}.
$$
Similarly we have the notation $\Delta_{<a}$ or $\Delta_{\geq a}$. We study geometric stability conditions in the parabolic region $\Delta_{<0}$. Denote by $L_{PE}$ the line passing through the two points $P$ and $E$ in the  $\cccp$. Recall a lemma due to Bayer \cite[Lemma 3]{LZ:Poisson}.

\begin{lemma} \label{lem:Bayer}
Let $P$ and $Q$ be two points in the region $\Delta_{<0}$ in the  $\cccp$. Let $F$ be a $\sigma_P$-stable object in $\coh_P$ with $\ch(F)\neq (0,0,1)$. Let $C$ and $D$ be the intersection points
$$ \{C, D\} \coloneqq L_{PF}\cap \Delta_0,
$$
of the line $L_{PF}$ and the parabola $\Delta_0$.
Denote the $\sigma_Q$-HN semistable factors of $F$ by $F_i$. Then for each factor, the phase $\phi_Q(F_i)$ lies in between $\phi_Q(C)$ and $\phi_Q(D)$.
\end{lemma}
\begin{proof}
The case that $\ch_0(F)\neq 0$ is proved in \cite[Lemma 3]{LZ:Poisson}.

So we assume that $\ch(F)=(0,\ch_1(F),\ch_2(F))$ with $\ch_1(F).H>0$. Now the point $F$ is in the $\infty$-line outside the $\cccp$. But the line $L_{PF}$ still makes sense: it is the line passing through the point $P$ with slope $\frac{\ch_2(F)}{\ch_1(F).H}$ in the $\cccp$, see \cite[Corollary 2.8]{Liu18}. So we still have the notation $l_{PF}^+$, which is the ray starting at the point $P$ on the line $L_{PF}$ with $s\geq s(P)$.
Note that $L_{QF}$ is parallel to $L_{PF}$. Then the proof follows by Li-Zhao's original argument.
\end{proof}

For a point $P=(1,s,q)$ in the  $\cccp$, we say that we move it along the parabola to the left by a number $b$ if we move it along the unique parabola of the form $\Delta_{a}$ passing through $P$ (so $a=\frac{1}{2}s^2-q$), and the result point is still on the same parabola $\Delta_a$ with $\frac{\ch_1}{\ch_0}$-coordinate $s-b$. Let $K$ be the canonical divisor of ${\PP}^2$, and $\omega_{{\PP}^2}$ be the dualizing sheaf.
Let $\sigma=\sigma_{s,q}$ with $(1,s,q)\in\Geo_{\mathrm{LP}}$. We identify $\sigma$ with the point $(1,s,q)$ in the  $\cccp$. Then $\sigma(-3)\coloneqq \sigma\otimes \omega_{{\PP}^2}$ is the point of moving $\sigma$ along the parabola to the left by $-H.K=3$. Similarly for $F\in \D$ as a point in the  $\cccp$, if $\ch_0(F)\neq 0$, then $F(-3)\coloneqq F\otimes \omega_{{\PP}^2}$ is the point of moving the point $F$ along the parabola to the left by $3$.

Let $A$, $B$, $\tilde{A}$, $\tilde{B}$ be the corresponding intersection points
$$
\{A, B\}\coloneqq L_{F\sigma}\cap\Delta_0, \quad \{\tilde{A},\tilde{B}\}\coloneqq  L_{F(-3)\sigma(-3)}\cap\Delta_0,
$$
with $s(B)>s(A)$ and $s(\tilde{B})>s(\tilde{A})$. We have the following observation.
\begin{lemma}\label{lem:horizental distance}
\begin{equation}\label{eq:same horizental difference}
s(B)-s(A) = s(\tilde{B})-s(\tilde{A}).
\end{equation}
\end{lemma}
\begin{proof}
This is an elementary calculation.
\end{proof}

We prove a lemma, which is the key calculation for proving $\gldim \sigma_{s,q} = 2$ in the region $\Delta_{<0}$.

\begin{lemma}\label{lem:above}
Let $\sigma_{s,q}$ be a geometric stability condition in the region $\Delta_{<0}$ in the $\cccp$. Denote $\sigma_{s,q}$ by $\sigma$.
Let $F,G$ be two $\sigma$ stable objects in a same heart satisfying:
$0<\phi_{\sigma}(F)<\phi_{\sigma}(G)\leq 1$, $\ch_0(F)\neq 0$ and $s<s(F)$.
Then $\Hom(F,G[2])=0$.
\end{lemma}
\begin{proof}
Let $P\coloneqq L_{F\sigma}\cap L_{F(-3)\sigma(-3)}$. We have two cases.

{\parindent =0pt \textbf{Case A.}}
$P$ is in the region $\Delta_{\geq 0}$.
Then by Lemma~\ref{lem:horizental distance},
we must have $s(\tilde{B})\leq s(A)$ (i.e. left above of Figure~\ref{fig:precise})
instead of $s(\tilde{B})>s(A)$ (i.e. left below or right below of Figure~\ref{fig:precise}).
So $l_{\sigma F}^+$ is above or equal to $l_{\sigma\tilde{A}}^+$ and $l_{\sigma\tilde{B}}^+$.

By \cite[Lemma A.3]{LLM19}, $F(-3)$ is $\sigma(-3)$-stable. By Lemma~\ref{lem:Bayer}, the $\sigma$-HN factor $F(-3)_i$ of $F(-3)$ lies between $\phi_{\sigma}(\tilde{A})$ and $\phi_{\sigma}(\tilde{B})$.
By \cite[Lemma 2]{LZ:Poisson}, $\phi_{\sigma}(F(-3)_i)\leq \phi_{\sigma}(F)$. So $\phi_{\sigma}^+(F(-3))\leq \phi_{\sigma}(F)<\phi_{\sigma}(G)$ and $\Hom(G,F(-3))=0$.
By Serre duality, we have
$$\Hom(F,G[2])\cong (\Hom(G[2],\Serre{(F)}))^*= (\Hom(G,F(-3)))^*=0.$$

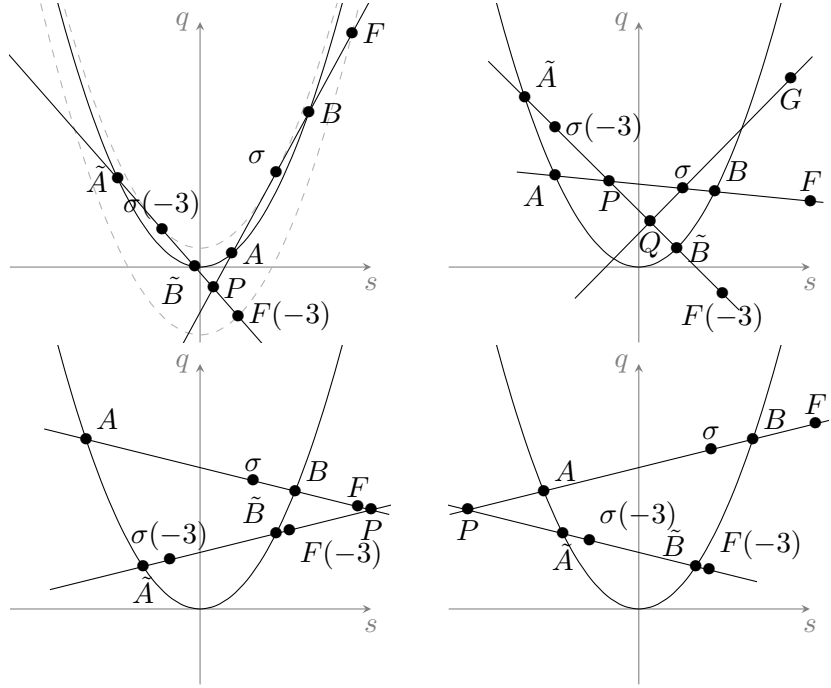
\begin{figure}
\centering
\begin{tikzpicture}[scale=.5,>=stealth]
\clip (-5,-2) rectangle (5,7);
  \draw[black!50, thin, ->] (0, -2) -- (0, 6.5) node[left]{$q$};
  \draw[black!50, thin, ->] (-5, 0) -- (4.5, 0) node[below]{$s$};
  \draw plot[smooth, domain=-5:5]   ( \x, .5*\x*\x);
  \draw[dashed,thin,black!30]   plot[smooth, domain=-5:5] ( \x, .5*\x*\x+.5);
  \draw[dashed,thin,black!30,]   plot[smooth, domain=-5:5] ( \x, .5*\x*\x-1.8);

  \draw (-1,1) coordinate ({Til_sigma});
  \draw (2,2.5) coordinate ({sigma});
  \draw (4,8-1.8) coordinate (F);
  \draw (1,.5-1.8) coordinate ({Til_F});
  \draw (2.861,4.0932) coordinate (B);
  \draw (0.8388,0.3518) coordinate (A);
  \draw (-0.1388,0.0096) coordinate ({Til_B});
  \draw (-2.1612,2.335) coordinate ({Til_A});

  \draw ($({sigma})!10!(F)$) --($(F)!10!({sigma})$);
  \draw ($({Til_sigma})!10!({Til_F})$) --($({Til_F})!10!({Til_sigma})$);

  \draw (intersection of  {sigma}--F and {Til_sigma}--{Til_F}) node {$\bullet$} node[right]{$P$};
  \draw({Til_sigma}) node{$\bullet$}node[above]{$\sigma(-3$)}
        ({sigma}) node{$\bullet$}node[above left]{$\sigma$}
        (F) node{$\bullet$}node[right]{$F$}
        ({Til_F}) node{$\bullet$}node[right]{$F(-3)$}
        ;
  \draw
        (A) node{$\bullet$}node[right]{${A}$}
        (B) node{$\bullet$}node[right]{${B}$}
        ({Til_A}) node{$\bullet$}node[left]{${\tilde{A}}$}
        ({Til_B}) node{$\bullet$}node[below left]{${\tilde{B}}$}
        ;
\end{tikzpicture}\qquad
\begin{tikzpicture}[yscale=0.5,xscale=0.5,>=stealth]
\clip (-5,-2) rectangle (5,7);
  \draw[black!50, thin, ->] (0, -2) -- (0, 6.5) node[left]{$q$};
  \draw[black!50, thin, ->] (-5, 0) -- (4.5, 0) node[below]{$s$};

  \draw  plot[smooth, domain=-5:5] ( \x, .5*\x*\x);

  \draw
    (-2.2,2.42) coordinate (A)
    (1,.5) coordinate ({Til_B})
    (-3,4.5) coordinate ({Til_A})
    (2,2) coordinate (B)
    ($(A)!.8!(B)$) coordinate ({sigma})
    ($(A)!1.6!(B)$) coordinate (F)
    ($({Til_A})!.2!({Til_B})$) coordinate ({Til_sigma})
    ($({Til_A})!1.3!({Til_B})$) coordinate ({Til_F})
    (4,5) coordinate (G);

  \draw ($({sigma})!1.1!(F)$) --($(F)!2.3!({sigma})$);
  \draw ($({Til_sigma})!1.1!({Til_F})$) --($({Til_F})!1.4!({Til_sigma})$);
  \draw ($({sigma})!1.2!(G)$) --($(G)!2!({sigma})$);

  \draw (intersection of  {sigma}--F and {Til_sigma}--{Til_F})node{$\bullet$}node[below]{$P$};

  \draw (intersection of  {sigma}--G and {Til_sigma}--{Til_F})node{$\bullet$}node[below]{${Q}$};

  \draw({Til_sigma}) node{$\bullet$}node[right]{$\sigma(-3$)}
        ({sigma}) node{$\bullet$}node[above]{$\sigma$}
        (F) node{$\bullet$}node[above]{$F$}
        ({Til_F}) node{$\bullet$}node[below]{$F(-3)$}
        (G) node{$\bullet$}node[below]{$G$}
        ;
  \draw
  (A) node{$\bullet$}node[below left]{$A$}
  (B) node{$\bullet$}node[above right]{$B$}
  ({Til_A}) node{$\bullet$}node[above right]{${\tilde{A}}$}
  ({Til_B}) node{$\bullet$}node[right]{${\tilde{B}}$};
\end{tikzpicture}\\
\begin{tikzpicture}[scale=.5,>=stealth]
\clip (-5,-2) rectangle (5,7);
  \draw[black!50, thin, ->] (0, -2) -- (0, 6.5) node[left]{$q$};
  \draw[black!50, thin, ->] (-5, 0) -- (4.5, 0) node[below]{$s$};

  \draw  plot[smooth, domain=-5:5] ( \x, .5*\x*\x);

  \draw(-3,4.5)coordinate (A);
  \draw(2.5,3.125)coordinate (B)
  ($(A)!.8!(B)$)coordinate({sigma})
  ($(A)!1.3!(B)$)coordinate(F);

  \draw(-1.5,1.125)coordinate (C);
  \draw(2,2)coordinate (D)
  ($(C)!.2!(D)$)coordinate(Q)
  ($(C)!1.1!(D)$)coordinate(G);

  \draw($({sigma})!1.3!(F)$) --($(F)!3!({sigma})$);
  \draw($(Q)!1.9!(G)$) --($(G)!2!(Q)$);

  \draw (intersection of  {sigma}--F and Q--G)node{$\bullet$}node[below]{$P$};

  \draw(Q) node{$\bullet$}node[above]{$\sigma(-3$)}
        ({sigma}) node{$\bullet$}node[above]{$\sigma$}
        (F) node{$\bullet$}node[above]{$F$}
        (G) node{$\bullet$}node[below right]{$F(-3)$};
  \draw
  (A) node{$\bullet$}node[above right]{$A$}
  (B) node{$\bullet$}node[above right]{$B$}
  (C) node{$\bullet$}node[below]{${\tilde{A}}$}
  (D) node{$\bullet$}node[above left]{${\tilde{B}}$};
\end{tikzpicture}\qquad
\begin{tikzpicture}[yscale=0.5,xscale=0.5,>=stealth]
\clip (-5,-2) rectangle (5,7);
 \draw[black!50, thin, ->] (0, -2) -- (0, 6.5) node[left]{$q$};
  \draw[black!50, thin, ->] (-5, 0) -- (4.5, 0) node[below]{$s$};

  \draw  plot[smooth, domain=-5:5] ( \x, .5*\x*\x);

  \draw(3,4.5)coordinate (B);
  \draw(-2.5,3.125)coordinate (A)
  ($(A)!.8!(B)$)coordinate({sigma})
  ($(A)!1.3!(B)$)coordinate(F);

  \draw(1.5,1.125)coordinate (D);
  \draw(-2,2)coordinate (C)
  ($(C)!.2!(D)$)coordinate(Q)
  ($(C)!1.1!(D)$)coordinate(G);

  \draw($({sigma})!1.2!(F)$) --($(F)!3.5!({sigma})$);
  \draw($(Q)!1.4!(G)$) --($(G)!2.2!(Q)$);

  \draw (intersection of  {sigma}--F and Q--G)node{$\bullet$}node[below]{$P$};

  \draw(Q) node{$\bullet$}node[above right]{$\sigma(-3$)}
        ({sigma}) node{$\bullet$}node[above]{$\sigma$}
        (F) node{$\bullet$}node[above]{$F$}
        (G) node{$\bullet$}node[above right]{$F(-3)$};
  \draw
  (A) node{$\bullet$}node[above right]{$A$}
  (B) node{$\bullet$}node[above right]{$B$}
  (C) node{$\bullet$}node[below]{${\tilde{A}}$}
  (D) node{$\bullet$}node[above left]{${\tilde{B}}$};
\end{tikzpicture}

\caption{Relative positions of $L_{F\sigma}$ and $L_{F(-3)\sigma(-3)}$: $P\in \Delta_{\geq 0}$ (left above); $P\in \Delta_{< 0}$ (right above). The below pictures are impossible by Lemma~\ref{lem:horizental distance}.}
\label{fig:precise}
\end{figure}

{\parindent =0pt \textbf{Case B.}}
$P$ is in the region $\Delta_{<0}$. So both $F$ and $F(-3)$ are $\sigma_P$-stable with
\begin{equation} \label{eq:phase-condition_R}
\phi_P(F)>\phi_P(F(-3)).
\end{equation}

Let $Q\coloneqq L_{G\sigma}\cap L_{F(-3)\sigma(-3)}$. We have three subcases.

{\parindent =0pt \textbf{Case B.(i)}}
$Q$ is in the region $\Delta_{>0}$. Then $Q$ is to the right of $\tilde{B}$ since $\tilde{B}$ is on the $\Delta_{0}$. Now $l_{QG}^+$ is above $l_{\sigma\tilde{A}}^+$ and $l_{\sigma\tilde{B}}^+$. We must have $l_{\sigma G}^+$ is above $l_{\sigma\tilde{A}}^+$ and $l_{\sigma\tilde{B}}^+$. By Lemma \ref{lem:Bayer} again, we have $\Hom(G,F(-3))=0$. By Serre duality, we have $\Hom(F,G[2])=0$.

{\parindent =0pt \textbf{Case B.(ii)}}
$Q$ is in the region $\Delta_{<0}$. We illustrate the picture in right above of Figure~\ref{fig:precise}.
Since $G$ is $\sigma$-stable, it is also $\sigma_{Q}$-stable. Since $F(-3)$ is $\sigma(-3)$-stable, it is also $\sigma_{Q}$-stable. We then compare their phases at $Q$ and have
$$
\phi_{Q}(G)=\phi_{\sigma}(G)>\phi_{\sigma}(F)=\phi_P(F)>\phi_P(F(-3))=\phi_{Q}(F(-3)),
$$
where each equality is because of colinear condition, and the first inequality is given by the assumption of the Lemma and the second inequality is given by (\ref{eq:phase-condition_R}).
So $\Hom(G,F(-3))=0$. By Serre duality, we have $\Hom(F,G[2])=0.$

{\parindent =0pt \textbf{Case B.(iii)}}
$Q$ is on the parabola $\Delta_{0}$. Since $F$ is $\sigma$-stable, we may perturb $\sigma$ a little bit and reduce to the previous cases.
\end{proof}

\begin{proposition}
Let $\sigma_{s,q}$ be in the region $\Delta_{<0}$ in the  $\cccp$. Then
\begin{equation}
\gldim \sigma_{s,q} = 2.
\end{equation}
\label{prop:gldiminParabola}
\end{proposition}
\begin{proof}
Denote $\sigma_{s,q}$ by $\sigma$. Let $F$ and $G$ be two $\sigma$-semistable objects such that 
$$\Hom(F, G[2])\neq 0.$$
Then by Serre duality,
\begin{equation}
\Hom(F, G[2]) \cong (\Hom (G, F(-3)))^*\neq 0.
\end{equation}
The object $F(-3)$ may not be $\sigma$-semistable. We consider its $\sigma$-HN factors.  Thus by \cite[Lemma 3.4]{Bri07} we have
$
\phi_{\sigma}(F)\leq \phi_{\sigma}(G[2])\leq \phi_{\sigma}^{+} (F(-3))+2
$.
So
\begin{equation}
0\leq \phi_{\sigma}(G[2]) - \phi_{\sigma}(F) \leq \phi_{\sigma}^{+}  (F(-3))-\phi_{\sigma}(F)+2.
\end{equation}
We need to show that
\begin{equation} \label{eq:phase-control}
\phi_{\sigma}(G[2]) - \phi_{\sigma}(F) \leq 2.
\end{equation}
The idea is to give an estimate of $\phi_{\sigma}^{+} (F(-3))-\phi_{\sigma}(F)$. We could assume that $F\in \coh_{\# s}$, i.e. $\phi_{\sigma}(F)\in (0,1]$.
We could also assume that $F$ is $\sigma$-stable since we can take its Jordan-H\"{o}lder factors. So by \cite[Lemma A.3]{LLM19}, $F(-3)$ is $\sigma(-3)$-stable.

We have the following three cases according to the Chern characters of $F$.

{\parindent =0pt \textbf{Case 1.}}
Assume $\ch_0(F)=0$, $\ch_1(F)=0$ and $\ch_2(F)>0$. Then $F$ is supported at point(s) and $\phi_{\sigma}(F(-3))=\phi_{\sigma}(F)$.
So (\ref{eq:phase-control}) holds.
On the other hand, for any closed point $x\in \PP^2$, we have $\Hom(\OO_x,\OO_x[2])\neq 0$ and
\begin{equation} \label{eq:gldim-reachable}
\phi_{\sigma}(\OO_x[2]) - \phi_{\sigma}(\OO_x)=2.
\end{equation}

{\parindent =0pt \textbf{Case 2.}}
Assume that $\ch_0(F)\neq 0$. We have the following three subcases.
\begin{itemize}
\item[(i)] $\sigma$ is to the left of $F$. This is precisely Lemma~\ref{lem:above}.

\item[(ii)]
If $\sigma$ is to the right of $F$, by applying a shifted derived dual functor, we reduce to case (i).

\item[(iii)]
If the $H$-slope of $F$ is $s$, by local finiteness of walls, we could replace $\sigma$ by $\sigma'$ in a small open neighbourhood of $\sigma$ so that $F$ is $\sigma'$-stable. So we reduce to case (i) or (ii).
\end{itemize}
{\parindent =0pt \textbf{Case 3.}}
Assume $\ch_0(F)=0$ and $\ch_1(F).H>0$. Now we have
$$\ch(F(-3))=(0,\ch_1(F),\ch_2(F)+\ch_1(F).K).$$

The line $L_{F\sigma}$ is the line passing through $\sigma$ of the slope $\frac{\ch_2(F)}{\ch_1(F).H}$.
Similarly, the line $L_{F(-3)\sigma(-3)}$ is the line passing through $\sigma(-3)$ of the slope $\frac{\ch_2(F)}{\ch_1(F).H}+H.K$ by \cite[Lemma A.3]{LLM19}. By Lemma~\ref{lem:Bayer}, the phase of $\phi_{\sigma}({(F(-3))}_i)$ lies between $\phi_{\sigma}(\tilde{A})$ and $\phi_{\sigma}(\tilde{B})$. We have similar analysis as the \textbf{Case 2} and still have  (\ref{eq:phase-control}).

Therefore for $\sigma\in\Delta_{<0}$ in the  $\cccp$, we have $\gldim(\sigma)=2$. Moreover, the value $2$ can be obtained by (\ref{eq:gldim-reachable}). This finishes the proof.
\end{proof}

%=========================================================
\section{Algebraic stability conditions} \label{sec:alg}
%=========================================================

%=========================================================
\subsection{Reviews} \label{sec:algReview}
We first recall the construction of algebraic stability
conditions with respect to exceptional triples from \cite{Li:spaceP2}.

\begin{definition}
We call an ordered set $\cE$ $=$ $\{E_1,E_2,E_3\}$
\emph{exceptional triple} on $\D^b(\PP^2)$ if $\cE$ is a
full strong exceptional collection of coherent sheaves on
$\D^b(\PP^2)$.
\label{def:exceptionaltriple}
\end{definition}

There is a one-to-one correspondence between the dyadic integers $\frac{p}{2^m}$ and exceptional bundles $E(\frac{p}{2^m})$:
\begin{equation*}
    \frac{p}{2^m} \Longleftrightarrow E(\frac{p}{2^m}),  \text{ for } p\in \ZZ \text{ and } m\in \ZZ_{\geq 0}.
\end{equation*}
The exceptional triples have been classified by Gorodentsev and Rudakov \cite{GorRu}. The exceptional triples are labeled by the following three cases,
$$
\Big\{\frac{p-1}{2^m},\frac{p}{2^m},\frac{p+1}{2^m}\Big\}, \quad \Big\{\frac{p}{2^m},\frac{p+1}{2^m},\frac{p-1}{2^m}+3\Big\}, \quad
\Big\{\frac{p+1}{2^m}-3,\frac{p-1}{2^m},\frac{p}{2^m}\Big\},
$$
for $p\in \ZZ$ and $m\in \ZZ_{\geq 0}$.
Note that the last two cases are mutations of the first case.

\begin{proposition} [{\cite[Section 3]{Mac07}}]
Let $\cE$ be an exceptional triple on $\D^b(\PP^2)$. For
any positive real numbers $m_1$, $m_2$, $m_3$ and real numbers
$\phi_1$, $\phi_2$, $\phi_3$ such that:
$$\phi_1<\phi_2<\phi_3,\text{ and }\phi_1+1<\phi_3,$$
there is a
unique stability condition $\sigma$ $=$ $(Z,\mathcal P)$ such that
\begin{enumerate}
\item each $E_j$ is stable with phase $\phi_j$;
\item $Z(E_j)=m_j e^{i\pi \phi_j}$.
\end{enumerate}
\label{prop:thetaE}
\end{proposition}

\begin{definition}
For an exceptional triple $\cE$ $=$ $\{E_1,E_2,E_3\}$ on
$\D^b(\PP^2)$, we write $\Theta_{\cE}$ as the space of all
stability conditions in Proposition \ref{prop:thetaE}, which is parametrized by
$$
S\coloneqq \{(m_1,m_2,m_3,\phi_1,\phi_2,\phi_3)\in (\mathbb R_{>0})^3\times
\mathbb R^3 \, | \,
\phi_1<\phi_2<\phi_3, \, \phi_1+1<\phi_3\}.
$$
We make the following
notations for some subsets of $\Theta_{\cE}$.
\begin{align*}
\Theta_{\cE}(A) &\coloneqq \{ \sigma \in \Theta_{\cE} | \sigma\in A\}, \text{ where } A \text{ is a subset of } S;\\
\Theta^{\mathrm{Pure}}_{\cE} &\coloneqq \{ \sigma \in \Theta_{\cE}|\phi_2-\phi_1\geq 1 \text{ and } \phi_3-\phi_2\geq1\};\\
\Theta^{\mathrm{left}}_{\cE, E_3} &\coloneqq \{ \sigma \in \Theta_{\cE}|\phi_2-\phi_1<1 \text{ and } E_3(3) \text{ is not } \sigma\text{-stable}\};\\
\Theta^{\mathrm{right}}_{\cE, E_1} &\coloneqq \{ \sigma \in \Theta_{\cE}|\phi_3-\phi_2<1 \text{ and } E_1(-3) \text{ is not } \sigma\text{-stable}\};\\
\Theta^{\Geo}_\cE &\coloneqq \Theta_\cE\cap \Stab^{\Geo}(\PP^2);\\
\Theta^{-}_{\cE, E_3} &\coloneqq \Theta_{\cE}(\phi_2-\phi_1<1) \setminus \Theta^{\Geo}_{\cE};\\
\Theta^{+}_{\cE, E_1} &\coloneqq \Theta_{\cE}(\phi_3-\phi_2<1) \setminus \Theta^{\Geo}_{\cE}.
\end{align*}
We denote
$$
\Stab^{\Alg}(\PP^2) \coloneqq \bigcup_{\cE \text{ exceptional triples}} \Theta_{\cE}
$$
and call the elements of it as the \emph{algebraic}
stability conditions. \label{def:thetaE}

\end{definition}

\begin{lemma}[{\cite[Lemma 2.4]{Li:spaceP2}}]
Let $\cE$ $=$ $\{E_1,E_2,E_3\}$ be an exceptional triple, and
$\sigma$ be a stability condition in $\Theta^{\mathrm{Pure}}_\cE$. The
only $\sigma$-stable objects are $E_i[n]$ for $i=1,2,3$ and $n\in
\mathbb Z$. \label{lemma:purepure}
\end{lemma}

%=========================================================
\subsection{Five points associated to an exceptional bundle}
%=========================================================
For an object $A\in\D$ with $\ch_0(A)\neq 0$, by abusing of notations, we write $A$ for $\vv(A)=(1,s(A),q(A))$ in (\ref{eq:s-q-coordinate}) as the associated point in the  $\cccp$, and call it the \emph{point} $A$ in the  $\cccp$. Moreover, by the Riemann--Roch formula, we have
\begin{equation}
    \chi(A,A)=\ch^2_0(A)(1-s(A)^2+2q(A)).
    \label{eq:chi_A_A}
\end{equation}
In particular, for an exceptional bundle $E$,  we have $\ch_0(E)\neq 0$, $\chi(E,E)=1$, and
\begin{equation} \label{eq:region of E}
   \frac{1}{2}s(E)^2-q(E) = \frac{1}{2} - \frac{1}{2\ch_0^2(E)}.
\end{equation}
So for each exceptional bundle $E$ as a point in the $\cccp$, the point $E$ is in the region $\Delta_{[0,\frac{1}{2})}$.

In the $\cccp$, for each exceptional bundle $E$, we define the following two pairs of parallel lines:
\begin{align}
    & \{ \chi(E,-)=0\},\quad
     \{ \chi(E,-)=\frac{\ch_0(-)}{\ch_0(E)}\};\\
    & \{ \chi(-,E)=0\},\quad
    \{ \chi(-,E)=\frac{\ch_0(-)}{\ch_0(E)}\}.
\end{align}

We now give a geometric description of above lines. By the Riemann--Roch formula, one can check that the line $\{ \chi(E,-)=\frac{\ch_0(-)}{\ch_0(E)}\}$ is the line $L_{E(-3)E}$ passing through the points $E(-3)$ and $E$. Similarly, the line $\{ \chi(-,E)=\frac{\ch_0(-)}{\ch_0(E)}\}$ is the line $L_{E(E(3))}$ passing through the points $E$ and $E(3)$.

The line $\{ \chi(E,-)=0\}$ is the line passing through points $E_1$ and $E_2$ for any choice of exceptional triple $\{E_1,E_2,E\}$ ending with $E$. It is clearly that this line is independent of the choice of $E_1$ and $E_2$. Similarly, the line $\{ \chi(-,E)=0\}$ is the line passing through points $E_2$ and $E_3$ for any choice of exceptional triple $\{E,E_2,E_3\}$ starting with $E$. This line is independent of the choice of $E_2$ and $E_3$.

For each exceptional bundle $E$, we define five points in the $\cccp$ as intersection points of the following lines or curves,
\begin{align*}
E^l &\coloneqq L_{E(E(3))}\cap\{\chi(E,-)=0\}, \\
E^r & \coloneqq L_{E(-3)E}\cap\{\chi(-,E)=0\},\\
E^+ &\coloneqq \{\chi(E,-)=0\}\cap\{\chi(-,E)=0\},\\
e^l & \coloneqq \Delta_{\frac{1}{2}}\cap \{\chi(E,-)=0\} \text{ as the first intersection point staring from } E^+,\\
e^r & \coloneqq \Delta_{\frac{1}{2}}\cap \{\chi(-,E)=0\} \text{ as the first intersection point staring from } E^+.  
\end{align*}

We now give a geometric description of above points. One can also refer to Figure~\ref{fig:regionMZc} and Figure~\ref{fig:MZvsM}. By the Riemann--Roch formula, (\ref{eq:region of E}) and (\ref{eq:s-q-coordinate}), we have
\begin{equation}
\label{eq:coordinate of E plus}
s(E^+)=s(E),\quad q(E^+)=q(E)-\frac{1}{(\ch_0(E))^2}.
\end{equation}
So $E^+$ is the point of moving $E$ downward of length $\frac{1}{(\ch_0(E))^2}$. By (\ref{eq:region of E}) and (\ref{eq:coordinate of E plus}), we have
\begin{equation} \label{eq:region of E plus}
   \frac{1}{2}s(E^+)^2-q(E^+) = \frac{1}{2} + \frac{1}{2\ch_0^2(E)}.
\end{equation}
So the point $E^+$ is in the region $\Delta_{(\frac{1}{2},1]}$. 

We observe that the point $E^l$ stands for $\vv(\mathsf{L}_{E}(E(3)))$, i.e. the reduced character of $\mathsf{L}_{E}(E(3))$. This is because by the definition (\ref{eq:leftmutation}), the point $\mathsf{L}_{E}(E(3))$ is on the line $L_{EE(3)}$. Also, the object $\mathsf{L}_{E}(E(3))\in E^\perp =\langle E_1,E_2\rangle$ has a resolution (\ref{eq:resolution}) (by taking $E_3=E$). Thus the point $\mathsf{L}_{E}(E(3))$ is on the line $\{\chi(E,-)=0\}$. 

By the Riemann--Roch formula, we have
\begin{equation*}
    \chi(E,E(3)) = 1+9\ch_0^2(E), \quad \chi(E(3),E)=1, \quad \chi(E(3),E(3))=1.
\end{equation*}
Since $[\mathsf{L}_{E}(E(3))]=[E(3)]-\chi(E,E(3))[E]$ in $K_{\mathrm{num}}(\PP^2)$, we have
\begin{equation*}
        \chi(\mathsf{L}_{E}(E(3)),\mathsf{L}_{E}(E(3))) = 1-\chi(E(3),E)\chi(E,E(3))  = -9\ch_0^2(E)<0.
\end{equation*}
Then by (\ref{eq:chi_A_A}), the point $E^l$ is in the region ${\Delta}_{>\frac{1}2}$. In particular, $E^l$ is in the line segment $\overline{E^+e^l}$. Similarly, ${E^r}$ stands for the reduced character of $\mathsf{R}_{E}(E(-3))$. It is in the region ${\Delta}_{>\frac{1}2}$ and in the line segment $\overline{E^+e^r}$. One can check that both of points $E^l$ and $E^r$ are in the parabola $\frac{1}{2}s^2-q = \frac{1}{2}+\frac{1}{18\ch_0^4(E)}$. 

\begin{definition}
\label{def:LPcurve}
(\cite[Definition 1.4]{Li:spaceP2}) The Le Potier curve $\mathrm{C}_{\mathrm{LP}}$ is a fractal curve defined in the  $\cccp$ as
\begin{equation*}
\mathrm{C}_{\mathrm{LP}} \coloneqq \bigsqcup_{\{E=E(\frac{p}{2^m}) \, | \, p\in \ZZ,\, m\in \ZZ_{\geq 0}\}} \left(\overline{E^+e^l}\cup \overline{E^+e^r}\right)\bigsqcup\{\text{Cantor pieces of } \Delta_{\frac{1}{2}}\}.
\end{equation*}
The region $\Geo_{\mathrm{LP}}$ is defined as $\Geo_{\mathrm{LP}}\coloneqq \Big\{(1,s,q)\in \cccp \, | \, (1,s,q)$ is above the curve $\mathrm{C}_{\mathrm{LP}}$ and is not on line segment $\overline{EE^+}$ for any exceptional bundle $E\Big\}$.
\end{definition}

%=========================================================
\subsection{Special regions associated to an exceptional triple}
%=========================================================

\begin{definition}
For an exceptional triple $\cE$ $=$ $\{E_1,E_2,E_3\}$, the region $\operatorname{MZ}^c_\cE$ is defined as the open region in the $\cccp$ bounded by the line segments $\overline{{E_1}{E_1^r}}$, $\overline{{E_1^r}{E_2}}$, $\overline{{E_2}{E_3^l}}$, $\overline{{E_3^l}{E_3}}$ and $\overline{{E_3}{E_1}}$ (see Figure~\ref{fig:regionMZc}).
The region $\operatorname{MZ}_{\cE}$ is defined as the open region in the $\cccp$ bounded by line segments $\overline{{E_1}{E_1^+}}$, $\overline{{E_1^+}{E_2}}$, $\overline{{E_2}{E_3^+}}$, $\overline{{E_3^+}{E_3}}$ and $\overline{{E_3}{E_1}}$ (see Figure~\ref{fig:MZvsM}). 

\label{def:me}
\end{definition}

\begin{figure}[ht]\centering
\begin{tikzpicture}[yscale=.75,xscale=.85,>=stealth]
\clip (-5,-2) rectangle (6,5);
\tikzset{%
    add/.style args={#1 and #2}{
        to path={%
 ($(\tikztostart)!-#1!(\tikztotarget)$)--($(\tikztotarget)!-#2!(\tikztostart)$)%
  \tikztonodes},add/.default={.2 and .2}}
}

\coordinate (E1) at (-1,0);
\coordinate (E2) at (0,-0.7);
\coordinate (E3) at (1,0.2);
\coordinate (F1) at (-3,4.5);
\coordinate (F3) at (3,4.7);

\draw [add= -1 and 0.25] (F1) to (E1) coordinate (ER1) node[below left]{$E_1^r$};
\draw [add= -1 and 0.27] (F3) to (E3) coordinate (EL3) node[below right ]{$E_3^l$};

\draw[white,fill=cyan!10] (E1)--(ER1)--(E2)--(EL3)--(E3)--cycle;

\draw[dashed] (F3) -- (E3);
\draw[dashed] (F1) -- (E1);
\draw[dashed] (E2) -- (E3);
\draw[dashed] (E2) -- (E1);

\draw[dashed] (E1) -- (ER1);
\draw[dashed] (ER1) -- (E2);
\draw[dashed] (E2) -- (EL3);
\draw[dashed] (EL3) -- (E3);
\draw[dashed] (E3) -- (E1);

\draw (E1) node {$\bullet$} node [left] {$E_1$};
\draw (E2) node {$\bullet$} node [above] {$E_2$};
\draw (E3) node {$\bullet$} node [right] {$E_3$};
\draw (F1) node {$\bullet$} node [below left] {$E_1(-3)$};
\draw (F3) node {$\bullet$} node [below right] {$E_3(3)$};
\draw (ER1) node {$\bullet$};
\draw (EL3) node {$\bullet$};

\end{tikzpicture}
\caption{The region of $\operatorname{MZ}^c_\cE$ in the  $\cccp$.}
\label{fig:regionMZc}
\end{figure}

We have $\operatorname{MZ}_{\cE}\subset \Geo_{\mathrm{LP}}$ and (\cite[Proposition 2.5]{Li:spaceP2}) $$
\Theta^{\Geo}_{\cE}=\glt\cdot\{\sigma_{s,q}\in\Stab^{\Geo}(\PP^2)\, | \, (1,s,q)\in \operatorname{MZ}_{\cE}\}.
$$

\begin{remark} Let $\cE=\{E_1,E_2,E_3\}$ be an exceptional triple. Note that the region $\operatorname{MZ}^c_\cE$ is a subregion of $\operatorname{MZ}_\cE$.
\begin{enumerate}
\item Since $E_2$ is in the region $\Delta_{[0,\frac{1}{2})}$ by (\ref{eq:region of E}) and $E_1^+, E_3^+$ are in the region $\Delta_{>\frac{1}{2}}$ by (\ref{eq:region of E plus}), we have
$$
e_1^r = \Delta_{\frac{1}{2}}\cap \overline{E_1^+E_2},\qquad e_3^l = \Delta_{\frac{1}{2}}\cap \overline{E_2E_3^+}.
$$

\item The line $L_{E_3(E_3(3))}$ is given as $\{\chi(-,E_3)=\frac{\ch_0(-)}{\ch_0(E_3)}\}$. For every stable vector bundle $A$ with slope between the slopes of $E_3$ and $E_3(3)$, we have $\chi(A,E_3)\leq 0$. The line segment $\overline{E_3E_3(3)}$ is contained in $\Geo_{\mathrm{LP}}$. \label{rem:meandmze1}

\item The point $E_3^l$ is on the line segment $\overline{e_3^lE_3^+}$. In particular, The reduced character of any exceptional bundles with slope smaller than that of $E_3$ is to the left of $e_3^l$.
\label{rem:meandmze2}

\item
By \cite[Corollary 1.19]{Li:spaceP2}, the exceptional object $E_3(3)$ is stable with respect to $\sigma_{s,q}$ for any $(1,s,q)$ in $\operatorname{MZ}^c_\cE$, and is destabilized by $E_3$ on the line segment $\overline{E_3E^l_3}$. In particular, the region $\operatorname{MZ}^c_\cE$ is a subregion of $\operatorname{MZ}_\cE$ by removing the region that either $E_3(3)$ or $E_1(-3)$ is not stable. In particular, we can identify the region $\operatorname{MZ}^c_\cE$ as the following algebraic stability conditions.
\label{rem:meandmze3}

\end{enumerate}

\label{rem:meandmze}
\end{remark}

\begin{lemma}
Let $\cE$ $=$ $\{E_1,E_2,E_3\}$ be an exceptional triple, then
$$\Theta_{\cE}\setminus (\Theta^{\mathrm{right}}_{\cE,E_1}\cup \Theta^{\mathrm{left}}_{\cE,E_3}\cup \Theta^{\mathrm{Pure}}_\cE) = \glt\cdot \Big\{\sigma_{s,q}\in\Stab^{\Geo}(\PP^2) \; | \; (1,s,q)\in
\operatorname{MZ}^c_\cE \Big\}.$$
\label{lemma:geoandalg}
\end{lemma}
\begin{proof}
By the previous Remark \ref{rem:meandmze}.\ref{rem:meandmze3}, the proof is the same as that for \cite[Lemma 1.29]{LZ:mmpP2}.
\end{proof}

\begin{definition}
For an exceptional triple $\cE$ $=$ $\{E_1,E_2,E_3\}$ on
$\D^b(\PP^2)$, we define $\operatorname{MZ}^l_{E_3}$ and $\operatorname{MZ}^r_{E_1}$ as subregions of $\operatorname{MZ}_{\cE}$ as follows:
\begin{align*}
\operatorname{MZ}^l_{E_3} &\coloneqq  \Big\{(1,s,q)\in \operatorname{MZ}_{\cE}\, |\, s<s(E_3), (1,s,q) \text{ is not above the line segment }\overline{E_3E^l_3}\Big\},\\
\operatorname{MZ}^r_{E_1} &\coloneqq  \Big\{(1,s,q)\in \operatorname{MZ}_{\cE}\, |\, s>s(E_1), (1,s,q) \text{ is not above the line segment }\overline{E_1E^r_1}\Big\}.
\end{align*}
\end{definition}{}

\begin{lemma}[Definition of $\Theta^{\mathrm{left}}_E$ and $\Theta^{\mathrm{right}}_E$] 
For any two exceptional triples $\cE$ and $\cE'$ on
$\D^b(\PP^2)$ ending with the same $E_3$ $=$ $E'_3$ $=$ $E$, we have
$\Theta^{\mathrm{left}}_{\cE,E_3}$ $=$ $\Theta^{\mathrm{left}}_{\cE',E'_3}$.
We denote \emph{this subspace} by $\Theta^{\mathrm{left}}_{E}$. In a similar way, we define
the subspace $\Theta^{\mathrm{right}}_E=\Theta^{\mathrm{right}}_{E_1}\coloneqq \Theta^{\mathrm{right}}_{\cE,E_1}$ for any exceptional triple $\cE$ starting with $E_1=E$. Moreover, we have
\begin{align}
\Theta^{\mathrm{left}}_{E_3}&= \Theta^{-}_{E_3} \bigsqcup \glt\cdot \Big\{\sigma_{s,q}\in\Stab^{\Geo}(\PP^2) \; | \; (1,s,q)\in \operatorname{MZ}^l_{E_3}\Big\},\label{eq:theta_E3_l}\\
\Theta^{\mathrm{right}}_{E_1}&= \Theta^{+}_{E_1} \bigsqcup \glt\cdot \Big\{\sigma_{s,q}\in\Stab^{\Geo}(\PP^2) \; | \; (1,s,q)\in \operatorname{MZ}^r_{E_1}\Big\}. \label{eq:theta_E1_r}
\end{align}
\label{prop:defoflegs}
\end{lemma}
\begin{proof}
By Remark \ref{rem:meandmze}.\ref{rem:meandmze3}, Lemma \ref{lemma:geoandalg} and \cite[Proposition and Definition 3.1]{Li:spaceP2}, we have the equation (\ref{eq:theta_E3_l}),
where $\Theta^{-}_{E_3}=\Theta^{-}_{\cE, E_3}$ is independent of the choice of $E_1$ and $E_2$. Note that by Remark \ref{rem:meandmze}.\ref{rem:meandmze2}, the boundary segment $\overline{E_3E^l_3}$ of $\operatorname{MZ}^l_{E_3}$ is also independent of the choice of $E_1$ and $E_2$ in the exceptional triple. The subspace $\Theta^{\mathrm{left}}_{E}$ is well-defined. Similarly, we have the equation (\ref{eq:theta_E1_r}).
\end{proof}

\begin{remark}
We illustrate the regions in Figure~\ref{fig:MZvsM}.
\begin{figure}[ht]\centering
\begin{tikzpicture}[yscale=.8,xscale=1,>=stealth]

\coordinate (E1) at (-3,0);
\coordinate (E2) at (0,-2.1);
\coordinate (E3) at (3,0.6);
\coordinate (E4) at (-3,-4.8);% E4=E_1^+
\coordinate (E5) at (3,-4.2); % E5=E_3^+
\coordinate (ER1) at (-1.245,-3.2205);
\coordinate (EL3) at (1.105,-2.874);

\draw[white,fill=green!10] (E3)--(EL3)--(E5)--cycle;
\draw[white,fill=green!10] (E1)--(ER1)--(E4)--cycle;
\draw[white,fill=cyan!10] (E1)--(ER1)--(E2)--(EL3)--(E3)--cycle;

\draw [red,dashed,domain=-1.2:1.1] plot (\x, {0.07*(\x+1.9)*(\x+1.9)-3});
\draw [red] (0,-3.3) node {$\Delta_{\frac{1}{2}}$};

\coordinate (eR1) at (-0.9263,-2.93367);
\coordinate (eL3) at (0.6417,-2.54919);

\coordinate (M) at (0,-0.7);
\draw [blue] (M) node {$\operatorname{MZ}^c_\cE$};

\coordinate (A) at (2.3,-2.5);
\draw [green] (A) node {$\operatorname{MZ}^l_{E_3}$};

\coordinate (B) at (-2.2,-2.5);
\draw [green] (B) node  {$\operatorname{MZ}^r_{E_1}$};

\draw[red,very thick] (eR1) -- (E4);
\draw[red,very thick] (eL3) -- (E5);

\draw[dashed] (E1) -- (E3);
\draw[green] (E3) -- (EL3);
\draw[dashed] (EL3) -- (E2);
\draw[dashed] (E2) -- (ER1);
\draw[green] (ER1) -- (E1);

\draw[dashed] (E1) -- (E4);
\draw[dashed] (E4) -- (ER1);

\draw[dashed] (E3) -- (E5);
\draw[dashed] (E5) -- (EL3);

\draw[dashed] (E2) -- (E3);
\draw[dashed] (E2) -- (E1);

\draw (E1) node {$\bullet$} node [below left] {$E_1$};
\draw (E2) node {$\bullet$} node [above] {$E_2$};
\draw (E3) node {$\bullet$} node [right] {$E_3$};
\draw (E4) node {$\bullet$} node [left] {$E_1^+$};
\draw (E5) node {$\bullet$} node [right] {$E_3^+$};

\draw (ER1) node {$\bullet$} node [below] {$E_1^r$};
\draw (EL3) node {$\bullet$} node [below] {$E_3^l$};

\draw [red] (eR1) node {$\bullet$} node[above] {$e_1^r$};
\draw [red] (eL3) node {$\bullet$} node[above] {$e_3^l$};

\end{tikzpicture}
\caption{The regions of $\operatorname{MZ}_\cE$, $\operatorname{MZ}^c_\cE$, $\operatorname{MZ}^l_{E_3}$ and $\operatorname{MZ}^r_{E_1}$ with relation (\ref{eq:regions MZ}) in the $\cccp$. The line segments $\overline{E_1^+e_1^r}$ and $\overline{E_3^+e_3^l}$ give parts of the Le Potier curve $\mathrm{C}_{\mathrm{LP}}$, and $\operatorname{MZ}_{\cE}\subset \Geo_{\mathrm{LP}}$.} 
\label{fig:MZvsM}
\end{figure}
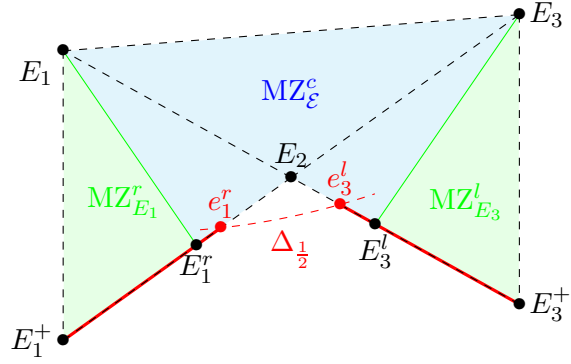
Then we could state Remark \ref{rem:meandmze}.\ref{rem:meandmze3} in a precise way, namely for an exceptional triple $\cE=\{ E_1, E_2, E_3 \}$,
\begin{equation}
  \operatorname{MZ}_\cE = \operatorname{MZ}^r_{E_1} \bigsqcup \operatorname{MZ}^c_\cE \bigsqcup \operatorname{MZ}^l_{E_3}.  \label{eq:regions MZ}
\end{equation}
\label{rem:MZvsM}
\end{remark}

%=========================================================
\section{Calculation of global dimension functions}
%=========================================================
The main result of this section is to compute the global dimension function on the algebraic stability conditions.
\begin{proposition}
Let $\cE=\{E_1,E_2,E_3\}$ be an exceptional triple on $\D^b(\PP^2)$ and $\Theta_\cE$ be the algebraic stability conditions with respect to $\cE$. The value of the global dimension function is

\begin{align*}
    \gldim(\sigma) = 
    \begin{cases}
    2, & \text{when }\sigma\in \Theta_{\cE}\setminus \left(\Theta^{\mathrm{right}}_{E_1}\cup \Theta^{\mathrm{left}}_{E_3}\cup \Theta^{\mathrm{Pure}}_\cE\right);\\
    \phi(\mathsf{R}_{E_1}(\Serre E_1))-\phi_1, &\text{when } \sigma \in \Theta^{\mathrm{right}}_{E_1};\\
    \phi_3-
    \phi(\mathsf{L}_{E_3}(\Serre^{-1} E_3)), &\text{when } \sigma \in \Theta^{\mathrm{left}}_{E_3};\\
    \phi_3-\phi_1, & \text{when } \sigma \in \Theta^{\mathrm{Pure}}_\cE.
    \end{cases}
\end{align*}
\label{prop:glonthestab}
\end{proposition}
Recall that $\mathsf{R}$ and $\mathsf{L}$ are the right and left mutations
in Section~\ref{sec:func}.
The rest of the section is devoted to the proof of the proposition above.

%=========================================================
\subsection{The locus with minimum global dimension}
%=========================================================
The other three cases are much more subtle, we first discuss the case when $\sigma\in \Theta_{\cE}\setminus \left(\Theta^{\mathrm{right}}_{E_1}\cup \Theta^{\mathrm{left}}_{E_3}\cup \Theta^{\mathrm{Pure}}_\cE\right)$.

\begin{proposition}

Let $\sigma$ be a stability condition in $\Theta_{\cE}\setminus \left(\Theta^{\mathrm{right}}_{E_1}\cup \Theta^{\mathrm{left}}_{E_3}\cup \Theta^{\mathrm{Pure}}_\cE\right)$, then $\gldim(\sigma)=2$.
\label{prop:mingldimpart}
\end{proposition}
The non-trivial part is the `$\leq$' part. As for a brief idea of the proof,  we will view $\sigma$ both as a stability condition in the region $\operatorname{MZ}^c_\cE$ in the  $\cccp$, and as a quiver stability condition. We will show that we only need to concern about $\Hom(F,G[2])\neq 0$ for two $\sigma$ stable objects $F$ and $G$ in a same heart with $\phi_{\sigma}(F)<\phi_{\sigma}(G)$. The line segments in the  $\cccp$ where $F$ and $G$ are $\sigma$ stable (i.e. $W_{F\sigma}$ and $W_{G\sigma}$ below) are `long' enough so that the line segment where $F(-3)$ is stable with respect to $\sigma(-3)$(i.e. $W_{F\sigma}(-3)$ below) intersects with previous two line segments $W_{F\sigma}$ and $W_{G\sigma}$. Then by the argument as that for stability conditions $\sigma_{s,q}$ above the parabola we show that $\Hom(F,G[2])=0$ and get a contradiction. Details of the proof is given as follows.

\begin{proof}[Proof for Proposition \ref{prop:mingldimpart}]
By Lemma \ref{lemma:geoandalg}, skyscraper sheaves are all stable with respect to $\sigma$. For any closed point $x\in \PP^2$, since $\Hom(\cO_x,\cO_x[2])=\mathbb{C}$, we have gldim$(\sigma)\geq 2$.

By Lemma \ref{lemma:geoandalg}, up to a $\glt$-action, we can view $\sigma$ as a stability $\sigma_{s,q}$ condition in the region $\operatorname{MZ}^c_\cE$ in the  $\cccp$. On the other hand, up to a  suitable $\mathbb C$-action on $\sigma$, we may let  the heart contain $E_1[2]$, $E_2[1]$ and $E_3$. Denote this stability condition and its heart by $\tilde{\sigma}$ and $\tilde{\cA}$ respectively.

\textbf{Step 1:} We reduce the equation in the proposition to the statement that for all stable objects $F$ and $G$ in $\tilde{\cA}$ with $\phi(F)<\phi(G)$, one must have $\Hom(F,G[2])= 0$.

 As $\{E_1[2],E_2[1],E_3\}$ is an Ext-exceptional collection (\cite[Definition 3.10]{Mac07}), an object in the heart is always of the form
$$E_1^{\oplus a_1}\rightarrow E_2^{\oplus a_2}\rightarrow E_3^{\oplus a_3}$$
for some non-negative integers $a_i$'s.

For any generators $E_i[3-i]$ in $\tilde{\cA}$, we always have
$$\Hom(E_i[3-i],E_j[3-j][m])=0$$
for every $m\geq 3$. Therefore, for any objects $F$ and $G$ in $\tilde{\cA}$, we have
$$\Hom(F,G[m])=0$$
for every $m\geq 3$. To prove the `$\leq$' part, we only need to show that for any $\sigma$-stable $F$ and $G$ with $\phi(F)<\phi(G)$ in the heart $\tilde{\cA}$, we have $\Hom(F,G[2])=0$.

\textbf{Step 2:} We show that the phases of $F$ and $G$ are both in $\left[ \phi(E_3(3)), \phi(E_1(-3)[2])\right]$.

 Suppose there are $\sigma$-stable $F$ and $G$ with $\phi(F)<\phi(G)$ in the heart $\tilde{\cA}$, such that  $\Hom(F,G[2])\neq 0$. Note that $\Hom(E_i[3-i],E_j[3-j][2])\neq 0$ if and only if $i=1$ and $j=3$, we must have $$\Hom(F,E_3[2])\neq 0\text{ and }\Hom(E_1[2],G[2])\neq 0.$$
By Serre duality, we have
$$\Hom(E_3(3),F)\neq 0\text{ and }\Hom(G,E_1(-3)[2])\neq 0.$$
By \cite[Corollary 1.19]{Li:spaceP2}, both objects $E_3(3)$ and $E_1(-3)[2]$ are  $\sigma_{s,q}$-stable (hence $\tilde{\sigma}$-stable). Both objects are in the heart $\tilde \cA$. Therefore, their phases satisfy the inequality:
\begin{equation}
\phi(E_3(3)) \leq \phi(F) < \phi(G) \leq \phi(E_1(-3)[2]).
    \label{eq:slopeofGFE}
\end{equation}

\textbf{Step 3:} We show that the walls $W_{F\sigma}$ and $W_{G\sigma}$ are `long' enough so that the wall $W_{F\sigma}(-3)$ intersects the walls $W_{F\sigma}$ and $W_{G\sigma}$. We compare their slopes and get the contradiction.

Here the \emph{wall} $W_{F\sigma}\coloneqq\{(1,s,q)\in \cccp|$ the line segment along the line $L_{F\sigma}$ that is above the Le Potier curve $\mathrm{C}_{\mathrm{LP}} \}$ and the wall
$$
W_{F\sigma}(-3)\coloneqq\{(1,s-3,q-3s+\frac{9}2)\; |\; (1,s,q)\in W_{F\sigma}\}.
$$
By Bertram's nested wall theorem, \cite[Corollary 1.24]{LZ:mmpP2}, the object $F$ is stable along the wall $W_{F\sigma}$. Let $F_a=(1, s(F_a), q(F_a))$ and $F_b=(1, s(F_b), q(F_b))$ be the two edges of the wall $W_{F\sigma}$ as that in the Figure~\ref{fig:compareslope1}. We denote similar notations for $G$ as that for $F$.

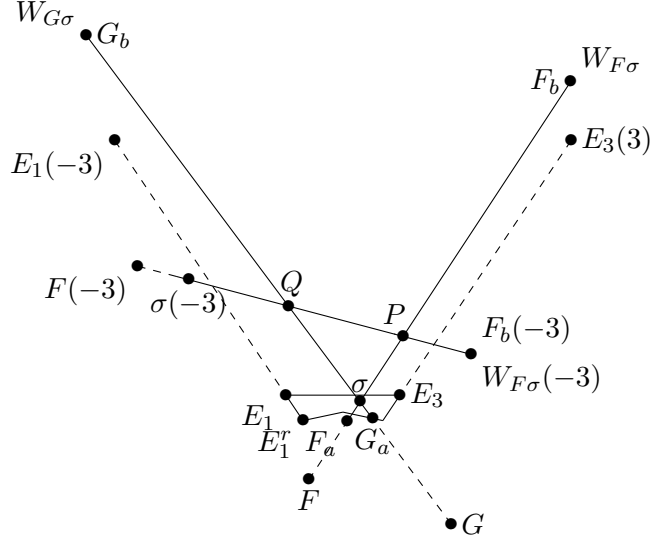
\begin{figure}[ht]\centering
\begin{tikzpicture}[yscale=.45,xscale=.75,>=stealth]
\clip (-6,-3.7) rectangle (7,12.2);
\tikzset{%
    add/.style args={#1 and #2}{
        to path={%
 ($(\tikztostart)!-#1!(\tikztotarget)$)--($(\tikztotarget)!-#2!(\tikztostart)$)%
  \tikztonodes},add/.default={.2 and .2}}
}
%position of sigma and F
\newcommand\xf{-0.6}
\newcommand\yf{-2}
\newcommand\xs{0.3}
\newcommand\ys{0.3}

\coordinate (E1) at (-1,0.5);
\coordinate (E2) at (0,0);
\coordinate (E3) at (1,0.5);
\coordinate (F1) at (-4,8); % E_1(-3)
%\coordinate (F2) at (-3,4.5);
%\draw (F2) node {$\bullet$} node [above left] {$E_2(-3)$};
%\coordinate (F33) at (-2,2);
%\draw (F33) node {$\bullet$} node [above left] {$E_3(-3)$};
\coordinate (F3) at (4,8); % E_3(3)
\coordinate (S) at (\xs,\ys); % sigma
\coordinate (SS) at (\xs-3,\ys-3*\xs+4.5); % sigma(-3)
\coordinate (G) at (1.9,-3.3);
\coordinate (F) at (\xf,\yf);
\coordinate (FF) at (\xf-3,\yf-3*\xf+4.5); % F(-3)
\coordinate (Q) at (intersection of  {SS}--{FF} and {G}--{S});
\coordinate (P) at (intersection of  {SS}--{FF} and {F}--{S});

\draw [add= -1 and 0.1] (F1) to (E1) coordinate (ER1) node[below left]{$E_1^r$};
\draw [add= -1 and 0.1] (F3) to (E3) coordinate (EL3);% node[below right ]{$E_3^l$};

\draw [add= 0 and -0.14, dashed] (G) to (S) coordinate (GS1) node {$\bullet$} node[below] {$G_a$};
\draw [add= -0.86 and 3] (G) to (S) coordinate (GS2) node {$\bullet$} node[right] {$G_b$} node[above left] {$W_{G\sigma}$};

\draw [add= 0 and -0.25, dashed] (FF) to (SS);
\draw [add= -0.75 and 5.5] (FF) to (SS) node {$\bullet$} node[above right] {$F_b(-3)$} node[below right] {$W_{F\sigma}(-3)$};

\draw [add= -5 and 4, dashed] (FF) to (SS) node[above] {$P$};
\draw [add= -3 and 2, dashed] (FF) to (SS) node[above] {$Q$};

\draw [add= 0 and -0.25, dashed] (F) to (S) coordinate (FS);
\draw [add= -0.75 and 4.1] (F) to (S) node {$\bullet$} node[left] {$F_b$} node[above right] {$W_{F\sigma}$};

\draw (FS) node {$\bullet$} node[below left] {$F_a$};

\draw (ER1) node {$\bullet$};
%\draw (EL3) node {$\bullet$};

\draw (E1) -- (E3);
%\draw (F1) -- (F33);
\draw[dashed] (F3) -- (E3);
\draw (ER1) -- (E2);
\draw (EL3) -- (E2);
\draw[dashed] (F1) -- (E1);

\draw (E1) node {$\bullet$} node [below left] {$E_1$};
%\draw (E2) node {$\bullet$} node [below] {$E_2$};
\draw (E3) node {$\bullet$} node [right] {$E_3$};
\draw (F1) node {$\bullet$} node [below left] {$E_1(-3)$};
\draw (F3) node {$\bullet$} node [right] {$E_3(3)$};
\draw (S) node {$\bullet$} node [above] {$\sigma$};
\draw (SS) node {$\bullet$} node [below] {$\sigma(-3)$};
\draw (G) node {$\bullet$} node [right] {$G$};
\draw (F) node {$\bullet$} node [below] {$F$};
\draw (FF) node {$\bullet$} node [below left] {$F(-3)$};
\draw (Q) node {$\bullet$};
\draw (P) node {$\bullet$};

\end{tikzpicture}
\caption{Compare the slopes of the wall $W_{G\sigma}$ and the wall $W_{F\sigma}(-3)$.}
\label{fig:compareslope1}
\end{figure}

By the relation of phases as (\ref{eq:slopeofGFE}), counter-clockwisely, one has the line segment $\overline{\sigma_{s,q}(E_3(3))}$, $\overline{\sigma_{s,q} F_b}$, $\overline{\sigma_{s,q} G_b}$ and $\overline{\sigma_{s,q}(E_1(-3))}$. In particular, either the wall $W_{F\sigma}$ is a vertical wall (parallel to the $\frac{\ch_2}{\ch_0}$-axis) or $|s(F_b)-s|>3$. Same statement holds for $W_{G\sigma}$. In every case, the segment $$\overline{\sigma_{s,q}(-3)F_b(-3)}=\overline{(1,s-3,q-3s+\frac{9}2)(1,s(F_b)-3,q(F_b)-3s(F_b)+\frac{9}2)}
$$
intersects both segments $\overline{\sigma_{s,q}F_b}$ and $\overline{\sigma_{s,q}G_b}$ at $P$ and $Q$ respectively. The object $F(-3)$ is stable at both $P$ and $Q$. By comparing the slopes, we have
$$\phi_Q(F(-3))=\phi_{P}(F(-3))<\phi_{P}(F)=\phi_{s,q}(F)<\phi_{s,q}(G)=\phi_Q(G).$$
By Serre duality,
$$\Hom(F,G[2])\cong(\Hom(G,F(-3)))^*=0.$$
We get the contradiction.
\end{proof}

%=========================================================
\subsection{The global dimension on the leg locus}
%=========================================================
We discuss the case that $\sigma\in \Theta^{\mathrm{left}}_{E_3}$. We first recall the following basic properties for an exceptional triple $\cE=\{E_1,E_2,E_3\}$.
Denote by $\mathrm{rk}(E)=\ch_0(E)$ and $\mathrm{hom}(E,F)=\dim \Hom(E,F)$.
\begin{lemma}
For an exceptional triple $\cE=\{E_1,E_2,E_3\}$, the ranks and homs of these exceptional objects satisfy the following equations.
\begin{align*}
    & (\mathrm{rk} E_1)^2+(\mathrm{rk} E_2)^2+(\mathrm{rk} E_3)^2=3 \mathrm{rk}E_1 \mathrm{rk} E_2 \mathrm{rk} E_3 \quad \text{(Markov equation)},\\
    & \mathrm{hom} (E_1,E_2)=3\mathrm{rk}E_3,\quad \mathrm{hom} (E_2,E_3)=3\mathrm{rk}E_1,\\
    & \mathrm{hom}(E_1,E_3)=9\mathrm{rk}E_1\mathrm{rk}E_3-3\mathrm{rk}E_2.
\end{align*}
The object $\mathsf{L}_{E_3}(E_3(3))[-1]$ admits a resolution:
\begin{equation}
    0\rightarrow E_1^{\oplus \mathrm{hom}(E_1,E_3)}\rightarrow E_2^{\oplus r}\rightarrow \mathsf{L}_{E_3}(E_3(3))[-1]\rightarrow 0, \label{eq:resolution}
\end{equation}
where $r=\mathrm{hom}(E_1,E_3)\mathrm{hom}(E_1,E_2)-\mathrm{hom}(E_2,E_3)$.
\label{lem:basicexc}
\end{lemma}
\begin{proof}
The equations of rank and hom are well-known in \cite{GorRu}. As for the last statement, we consider the resolution of $E_3(3)$. Note that $\cD_\infty(\PP^2)$ has the semiorthorgonal decomposition $\langle E_1,E_2,E_3\rangle$, so an object $A$ admits a unique filtration
$$0=F_0\subset F_1\subset F_2\subset F_3=A$$
such that $\Cone(F_{i}\rightarrow F_{i+1})\in \langle E_{3-i}\rangle $ for $i=0,1,2$. The term $\Cone(F_{0}\rightarrow F_1)$ is given by $\bigoplus_iE_3[i]\otimes \Hom(E_3[i],A)$, while the term $\Cone(F_{2}\rightarrow F_3)$ is given by $\bigoplus_iE_1[i]\otimes \Hom(A,E_1[i])^*$.

When $A=E_3(3)$, we have $\Cone(F_{0}\rightarrow F_1)= E_3\otimes \Hom(E_3,E_3(3))= E_3^{\oplus 9(\mathrm{rk}E_3)^2+1}$ and $\Cone(F_{2}\rightarrow F_3)= E_1^{\oplus \mathrm{hom}(E_1,E_3)}[2]$. The factor $\Cone(F_{1}\rightarrow F_2)$ can only be $E_2^{\oplus r}[1]$. By the equations of rank and hom in the lemma, the rank $$r=\frac{9(\mathrm{rk}E_3)^3+9\mathrm{rk}E_3(\mathrm{rk}E_1)^2}{\mathrm{rk} E_2}-3\mathrm{rk}E_1=\mathrm{hom}(E_1,E_3)\mathrm{hom}(E_1,E_2)-\mathrm{hom}(E_2,E_3).$$
Note that $\mathsf{L}_{E_3}(E_3(3))[-1]$ is the kernel of the map $E_3\otimes \Hom(E_3,E_3(3))\xrightarrow{\mathrm{ev}}E_3(3)$, the resolution sequence is clear.
\end{proof}

\begin{lemma}
Let $\sigma$ be a stability condition in $\Theta_\cE(\phi_2<\phi_1+1)$, suppose an object $F=\Cone(E_1^{\oplus a}\rightarrow E_2^{\oplus b})$ is stable with respect to $\sigma$, then $F$ is stable everywhere in $\Theta_\cE(\phi_2<\phi_1+1)$.
\label{lem:twotermstabobj}
\end{lemma}
\begin{proof}
For any stability condition in $\Theta_\cE(\phi_2<\phi_1+1)$, by a suitable $\mathbb C$-action, we may assume that the heart contains $E_1[2]$, $E_2[1]$ and $E_3[n]$ for some $n\leq 0$. As $\{E_1[2],E_2[1],E_3[n]\}$ is an Ext-exceptional collection, an object in the heart is always of the form
$$E_1^{\oplus a_1}\rightarrow E_2^{\oplus a_2}\rightarrow E_3^{\oplus a_3}.$$

The object $F[1]$ can only be destabilized  by some subobjects $F'=\Cone(E_1^{\oplus a'}\rightarrow E_2^{\oplus b'})[1]$ in the heart generated by $\{ E_1[2],E_2[1],E_3[n]\}$ with larger phase, which means $\frac{a'}{b'}>\frac{a}b$. Note that this is independent of the choice of $\sigma$ in $\Theta_\cE(\phi_2<\phi_1+1)$, the object $F$ is stable everywhere in $\Theta_\cE(\phi_2<\phi_1+1)$.
\end{proof}

Now we are ready to compute the example achieving the value of the global dimension function in the region of $\Theta^{\mathrm{left}}_{E_3}$.

\begin{lemma}
Let $\sigma$ be a stability condition in $\Theta^{\mathrm{left}}_{E_3}$, then $\mathsf{L}_{E_3}E_3(3)$ is $\sigma$-stable and it has a non-zero morphism to $E_3[2]$. In particular, we have $\gldim(\sigma)\geq \phi_3-
    \phi(\mathsf{L}_{E_3}E_3(3))+2$.
\label{lem:lowbdforgldonlegs}
\end{lemma}

\begin{proof} We are working in the $\cccp$.
By \cite[Corollary 3.2]{LZ:mmpP2}, the object  $$\mathsf{L}_{E_3}E_3(3)=\Cone(E_3\otimes \Hom(E_3,E_3(3))\xrightarrow{\mathrm{ev}} E_3(3))$$
is $\sigma_{s,q}$-stable for $(1,s,q)$ which is slightly above the line segment $\overline{E_3E_3(3)}$. The object $\mathsf{L}_{E_3}E_3(3)$ is stable along the line segment $\overline{(1,s,q)E_3^l}$. As this segment intersects $\operatorname{MZ}^c_\cE$, by Lemma  \ref{lemma:geoandalg}, the object $\mathsf{L}_{E_3}E_3(3)$ is stable with respect to some stability condition in $\Theta_\cE(\phi_2<\phi_1+1)$. By Lemma \ref{lem:basicexc} and \ref{lem:twotermstabobj}, the object $\mathsf{L}_{E_3}E_3(3)$ is $\sigma$-stable for every $\sigma\in \Theta_\cE(\phi_2<\phi_1+1)$.

By applying $\Hom(-,E_3[2])$ on the distinguished triangle
$$
E_3\otimes \Hom(E_3,E_3(3))\xrightarrow{\mathrm{ev}}E_3(3)\rightarrow \mathsf{L}_{E_3}E_3(3)\xrightarrow{+},
$$
we have
$\Hom(\mathsf{L}_{E_3}E_3(3),E_3[2])\cong\Hom(E_3(3),E_3[2])=\mathbb C$.
\end{proof}

As for the `$\leq$' direction, we first treat with the easier case that the stable objects can be classified.
\begin{proposition}
Let $\sigma$ be a stability condition in $\Theta_\cE(\phi_2<\phi_1+1< \phi_3-1)$, then up to a homological shift, a $\sigma$-stable object is either
\begin{itemize}
    \item $E_3$ or
    \item $\Cone(E_1^{\oplus a}\rightarrow E_2^{\oplus b})$
\end{itemize}
induced by a stable quiver representation $\mathbb{C}^{\oplus a}\xRrightarrow{\mathrm{hom}(E_1,E_2)\text{ arrows}} \mathbb{C}^{\oplus b}$. Moreover, we have $\gldim(\sigma)=\phi_3-\phi(\mathsf{L}_{E_3}E_3(3))+2$.
\label{prop:stabobjinfarleg}
\end{proposition}
\begin{proof}
 By a suitable $\mathbb C$-action, we may assume that the heart contains $E_1[2]$, $E_2[1]$ and $E_3[n]$ for some $n\leq -1$. As $\{E_1[2],E_2[1],E_3[n]\}$ is an Ext-exceptional collection, an object in the heart is $\Cone(E_1^{\oplus a}\rightarrow E_2^{\oplus b})[1]\bigoplus E_3^{\oplus c}[n]$. An object $\Cone(E_1^{\oplus a}\rightarrow E_2^{\oplus b})[1]$ is $\sigma$-stable if and only if for any non-zero proper subobject $\Cone(E_1^{\oplus a_1}\rightarrow E_2^{\oplus b_1})[1]$ we have $\frac{a_1}{b_1}<\frac{a}b$. The first part of the statement is clear.

As for the second part of the statement, by Lemma \ref{lem:lowbdforgldonlegs}, we only need to show the `$\leq$' side, note that for any two stable objects $F$ and $F'$ in the form of $\Cone(E_1^{\oplus a}\rightarrow E_2^{\oplus b})$, we always have $\Hom(F,F'[m])=0$ for $m\geq 2$. By the classification of stable objects, we only need to consider potential non-zero morphisms from $\Cone(E_1^{\oplus a}\rightarrow E_2^{\oplus b})$ to $E_3[m]$ for $m\geq 1$. When $\phi(\Cone(E_1^{\oplus a}\rightarrow E_2^{\oplus b}))<\phi(\mathsf{L}_{E_3}E_3(3)[-1])$, by Lemma \ref{lem:basicexc}, we have
\begin{equation}
    \frac{b}a>\mathrm{hom}(E_1,E_2)-\frac{\mathrm{hom}(E_2,E_3)}{\mathrm{hom}(E_1,E_3)}. \label{eq:sloptwoterms}
\end{equation}
Let $\phi_3'$ be $\phi(\mathsf{L}_{E_3}E_3(3))$, which is greater than $\phi_1+1$. We consider the stability condition $\sigma'$ in $\Theta_\cE$ given by $(m_1,m_2,m_3,\phi_1,\phi_2,\phi_3')$.  By Lemma \ref{lem:twotermstabobj}, $\Cone(E_1^{\oplus a}\rightarrow E_2^{\oplus b})$ is $\sigma'$-stable and $E(3)$ is $\sigma'$-semistable with phase $\phi_3'=\phi'(E_3)=\phi'(\mathsf{L}_{E_3}E_3(3))$. By (\ref{eq:sloptwoterms}), we have $\phi'(\Cone(E_1^{\oplus a}\rightarrow E_2^{\oplus b}))<\phi'(\mathsf{L}_{E_3}E_3(3)))-1=\phi'(E(3))-1$. Therefore, for any $m\geq 1$, by Serre duality, we have
$$\Hom(\Cone(E_1^{\oplus a}\rightarrow E_2^{\oplus b}),E_3[m])\cong (\Hom(\Cone(E_3(3)[m-2],E_1^{\oplus a}\rightarrow E_2^{\oplus b})))^*=0.$$
As a summary, the global dimension at $\sigma$ is $\phi_3-
    \phi(\mathsf{L}_{E_3}E_3(3))+2$, and is achieved via the morphism between $\mathsf{L}_{E_3}E_3(3)$ and $E_3[2]$.
\end{proof}

We finally treat with region $\Theta_\cE(\phi_2<\phi_3-1<\phi_1+1)\cap \Theta^{\mathrm{left}}_{E_3}$, where the stable objects are more complicated. In this case, the potential stable characters are away from the kernel of central charge of every $\sigma$ in $\Theta_\cE(\phi_2<\phi_3-1<\phi_1+1)$. We will think both the stable characters and (kernels of central charges of) stability conditions in the  $\cccp$. This will allow us to show the vanishing of certain morphisms by comparing slopes.

We first prove a nested wall result for the algebraic stability conditions. Denote $\Theta^+_\cE(\phi_2<\phi_1+1,\phi_3<\phi_1+2)\coloneqq\{\sigma\in \Theta_\cE(\phi_2<\phi_1+1,\phi_3<\phi_1+2)|$  the kernel of central charge of $\sigma$ is spanned by $(1,s,q)$ for some $s>s(E_1)\}$.

\begin{lemma}
Let $\sigma$ be a stability condition in $\Theta^+_\cE(\phi_2<\phi_1+1,\phi_3<\phi_1+2)$ and  $G$ be a $\sigma$-stable object. Then $G$ is $\sigma'$-stable for every $\sigma'$ in $\Theta^+_\cE(\phi_2<\phi_1+1,\phi_3<\phi_1+2)$ with kernel of central charge on the line through $G$ and $\sigma$.
\label{lem:nestedforalgstab}
\end{lemma}
\begin{proof}
In the $\cccp$, the kernel of the central charge of $\sigma$ is in the region bounded by rays through $E_1E_3$, $E_1E_2$ as shown in the Figure~\ref{fig:stabsigmaF} (Area I).

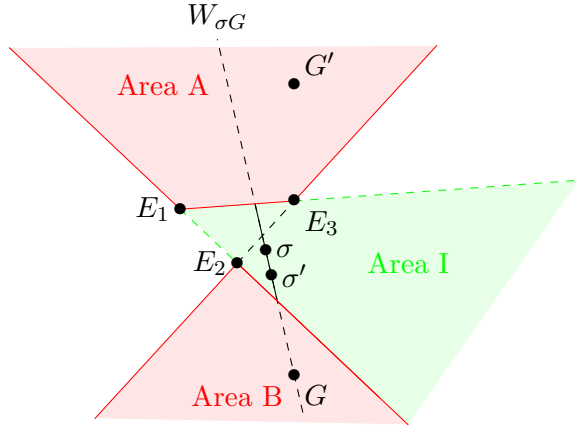
\begin{figure}[ht]\centering
\begin{tikzpicture}[yscale=.55,xscale=.75,>=stealth]
\clip (-6,-5.5) rectangle (6,7);
\tikzset{%
    add/.style args={#1 and #2}{
        to path={%
 ($(\tikztostart)!-#1!(\tikztotarget)$)--($(\tikztotarget)!-#2!(\tikztostart)$)%
  \tikztonodes},add/.default={.2 and .2}}
}

\coordinate (E1) at (-1,0);
\coordinate (E2) at (0,-1.3);
\coordinate (E3) at (1,0.2);

\coordinate (G) at (1,-4);
\coordinate (C) at (0.5,-1); % C=sigma

\draw[white,fill=red!10] ($(E1)!-3!(E2)$)--(E1)--(E3)--($(E3)!-2.5!(E2)$)--cycle;
\draw[white,fill=red!10] ($(E1)!4!(E2)$)--(E2)--($(E3)!3.5!(E2)$)--cycle;
\draw[white,fill=green!10] ($(E1)!3.5!(E3)$)--(E1)--(E2)--($(E1)!4!(E2)$)--cycle;

\draw (G) node {$\bullet$} node [below right] {$G$};
\draw (1,3) node {$\bullet$} node [above right] {$G'$};

\draw (C) node {$\bullet$} node [right] {$\sigma$};

\draw [add=0.3 and 1.7, dashed] (G) to (C) node[above] {$W_{\sigma G}$};
\draw [add=-0.61 and 0.37] (G) to (C);
\draw [add=-0.8 and -0.2, dashed] (G) to (C) node {$\bullet$} node[right] {$\sigma'$};

\draw[red] (E2) -- ($(E1)!4!(E2)$);
\draw[red] (E1) -- ($(E1)!-3!(E2)$);

\draw [add= -1 and 2.5, red] (E2) to (E3);
\draw [add= -1 and 2.5, red] (E1) to (E2);
\draw [add= -1 and 2.5, red] (E3) to (E2);
\draw [dashed,green] (E1) -- (E2);

\draw [red] (-1.3,3) node {Area A};
\draw [red] (0,-4.5) node {Area B};
\draw [green] (3,-1.3) node {Area I};

\draw [red] (E1) -- (E3);
\draw[dashed] (E2) -- (E3);
\draw[dashed,green] (E3) -- ($(E1)!3.5!(E3)$);

\draw (E1) node {$\bullet$} node [left] {$E_1$};
\draw (E2) node {$\bullet$} node [left] {$E_2$};
\draw (E3) node {$\bullet$} node [below right] {$E_3$};

\end{tikzpicture}
\caption{Stability conditions through $W_{\sigma G}$.}
\label{fig:stabsigmaF}
\end{figure}

By a suitable $\mathbb C$-action, we may assume that the heart contains $\{E_1[2],E_2[1],E_3\}$. Denote this heart by $\tilde{\cA}$, then an object in $\tilde{\cA}$ is of the form $E_1^{\oplus a_1}\rightarrow E_2^{\oplus a_2}\rightarrow E_3^{\oplus a_3}.$ In particular, the reduced character of a stable object is in the closed region (Area A $\cup$ Area B in Figure~\ref{fig:stabsigmaF}) bounded by the rays through $E_1E_2$, $E_2E_3$ and line segment through $E_1E_3$.

The phase of  $G$ is determined by the slope of line through $\sigma$ and $\vv(G)$. As for another object $G'$, its phase $\phi(G')<\phi(G)$ if and only if the line through $\sigma$ and $\vv(G')$ rotates counter-clockwisely to the line through $\sigma$ and $\vv(G)$ without passing though the line through $\sigma$ and $E_1[2]$.

For every non-zero proper subobject $G'$ of $G$ in $\tilde{\cA}$, since $G$ is stable, $G'$ has smaller phase than that of $G$. In the $\cccp$, that is equivalent to the following description for $\vv(G')$:

The reduced character of $G'$ is either to the right of the line through $G$ and $\sigma$ when it is in Area A, or it is to the left of the line through $G$ and $\sigma$ when it is in Area B.

Note that for every stability condition $\sigma'$ in $\Theta^+_\cE(\phi_2<\phi_1+1,\phi_3<\phi_1+2)$ with kernel of central charge on the line through $G$ and $\sigma$, the line through $G'$ and $\sigma'$ rotates counter-clockewisely to the line through $\sigma$, $\sigma'$ and $G$. The object $G$ is $\sigma'$-stable.
\end{proof}

\begin{proposition}
Let $\sigma$ be a stability condition in $\Theta_\cE(\phi_2<\phi_3-1<\phi_1+1)\cap \Theta^{\mathrm{left}}_{E_3}$. Then $\gldim(\sigma)=\phi(E_3)-
    \phi(\mathsf{L}_{E_3}E_3(3))+2$.
\label{prop:glfunconhipleg}
\end{proposition}
\begin{proof}
By Lemma \ref{lem:lowbdforgldonlegs}, we only need to show the `$\leq$' part. 

By a suitable $\mathbb C$-action, we may assume that the heart contains $\{E_1[2],E_2[1],E_3\}$. Denote this heart by $\tilde{\cA}$, we have the same description for objects in $\tilde{\cA}$ as that in Lemma \ref{lem:nestedforalgstab}.

\textbf{Step 1:} We reduce the claim in the proposition to the following statement: for all stable objects $F$ and $G$ in $\tilde{\cA}$ with $\Hom(F,G[2])\neq 0$, the difference of their phases $\phi(G)-\phi(F)\leq \phi(E_3)-\phi(\mathsf{L}_{E_3}E_3(3))$.

In the $\cccp$, the kernel of the central charge of $\sigma$ is in the region bounded by rays through $E_1E_3$, $E_1E_2$ and line segment $\overline{E_3E_3^l}$ as shown in the Figure~\ref{fig:regofquiv} (Area I $\cup$ Area II). Recall that the point $E_3^l$ and $\mathsf{L}_{E_3}(E_3(3))$ are the same point in the $\cccp$.

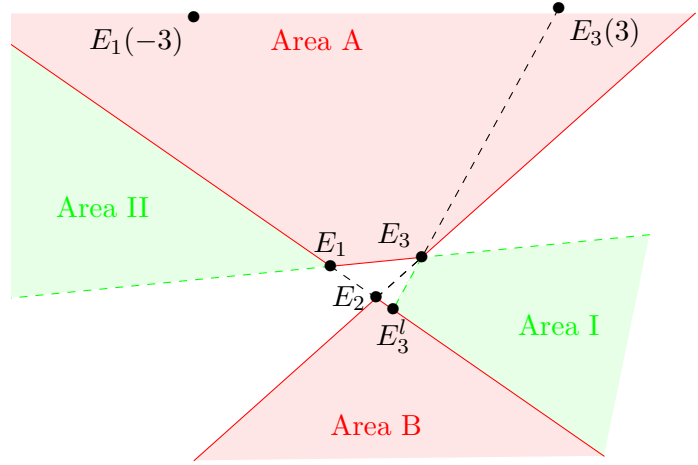
\begin{figure}[ht]\centering
\begin{tikzpicture}[yscale=.6,xscale=.6,>=stealth]
\clip (-8,-4.5) rectangle (8,6);
\tikzset{%
    add/.style args={#1 and #2}{
        to path={%
 ($(\tikztostart)!-#1!(\tikztotarget)$)--($(\tikztotarget)!-#2!(\tikztostart)$)%
  \tikztonodes},add/.default={.2 and .2}}
}

\coordinate (E1) at (-1,0);
\coordinate (E2) at (0,-0.7);
\coordinate (E3) at (1,0.2);
\coordinate (F1) at (-4,5.5);% F1=E_1(-3)
\coordinate (F3) at (4,5.7);% F3=E_3(3)

\draw [add= -1 and 0.21, green] (F3) to (E3) coordinate (EL3);

\draw[white,fill=red!10] ($(E1)!-8!(E2)$)--(E1)--(E3)--($(E3)!-6!(E2)$)--cycle;
\draw[white,fill=green!10] ($(E1)!-7!(E2)$)--(E1)--($(E3)!6!(E1)$)--cycle;
\draw[white,fill=red!10] ($(E1)!6!(E2)$)--(E2)--($(E3)!5!(E2)$)--cycle;
\draw[white,fill=green!10] ($(E1)!3.5!(E3)$)--(E3)--(EL3)--($(E1)!6!(E2)$)--cycle;

\draw[dashed] (F3) -- (E3);
\draw[dashed] (E2) -- (E3);
\draw[dashed] (E2) -- (E1);

\draw[red] ($(E1)!-8!(E2)$)--(E1);
\draw[red] (E1) -- (E3);
\draw[red] (E3)--($(E3)!-6!(E2)$);
\draw[red] ($(E1)!6!(E2)$)--(E2);
\draw[red] (E2)--($(E3)!5!(E2)$);

\draw[dashed,green] (E1)--($(E3)!6!(E1)$);

\draw[dashed,green] ($(E1)!3.5!(E3)$)--(E3);
\draw[dashed,green] (E3)--(EL3);

\draw [red] (-1.3,5) node {Area A};
\draw [red] (0,-3.5) node {Area B};
\draw [green] (4,-1.3) node {Area I};
\draw [green] (-6,1.3) node {Area II};

\draw[dashed] (E2) -- (E3);

\draw (E1) node {$\bullet$} node [above] {$E_1$};
\draw (E2) node {$\bullet$} node [left] {$E_2$};
\draw (E3) node {$\bullet$} node [above left] {$E_3$};
\draw (F1) node {$\bullet$} node [below left] {$E_1(-3)$};
\draw (F3) node {$\bullet$} node [below right] {$E_3(3)$};
\draw (EL3) node {$\bullet$} node[below] {$E_3^l$};

\end{tikzpicture}
\caption{Stable characters are in Area A $\cup$ Area B. The kernels of the central charges are in Area I $\cup$ Area II.}
\label{fig:regofquiv}
\end{figure}

As for any generators in $\{E_1[2],E_2[1],E_3\}$, we have $\Hom(-,-[m])=0$ for any $m\geq 3$. For any objects $F$ and $G$ in the heart, we have $\Hom(F,G[m])=0$ for any $m\geq 3$. To prove the `$\leq$' part of the statement, we only need consider $\Hom(F,G[2])\neq 0$ for stable objects $F,G$ in the heart with $\phi(F)<\phi(G)$.

Suppose there are $\sigma$-stable objects $F$ and $G$ with
\begin{equation}
\phi(G)-\phi(F)>\phi(E_3)-\phi(\mathsf{L}_{E_3}E_3(3))
    \label{eq:comparephFG}
\end{equation}
 in the heart $\tilde{\cA}$, such that $\Hom(F,G[2])\neq 0$. By the same argument as that in Proposition \ref{prop:mingldimpart}, we must have
\begin{equation}
    \Hom(E_3(3),F)\neq 0\text{ and }\Hom(G,E_1(-3)[2])\neq 0.
    \label{eq:nonvanofe3fge1}
\end{equation}

\textbf{Step 2:} We show that $\phi(G)>\phi(E_3)$.

Suppose $\phi(F)<\phi(\mathsf{L}_{E_3}E_3(3))$, then $\phi(F)<\phi(E_3)$. Therefore, the object $F$ is of the form $\Cone(E_1^{\oplus a_F}\rightarrow E_2^{\oplus b_F})[1]$. 
By Lemma \ref{lem:twotermstabobj}, $F$ is stable everywhere in $\Theta(\phi_1+1>\phi_2)$. In particular, it is stable with every stability condition $\sigma'$ on the line segment $\overline{E_3 E_3^l}$, where $E_3(3)$ is $\sigma'$-semistable. Since $\Hom(E_3(3),F)\neq 0$, we have $\phi'(F)\geq \phi'(E_3(3))=\phi'(\mathsf{L}_{E_3}E_3(3))$. 
Therefore, we have
$$ \frac{b_F}{a_F}\leq \mathrm{hom}(E_1,E_2)-\frac{\mathrm{hom}(E_2,E_3)}{\mathrm{hom}(E_1,E_3)}, \quad \phi'(F)\geq \phi'(\mathsf{L}_{E_3}E_3(3)),$$
which contradicts the assumption that $\phi(F)<\phi(\mathsf{L}_{E_3}E_3(3))$.

By (\ref{eq:comparephFG}), we must have
\begin{equation}
    \phi(G)>\phi(E_3).
    \label{eq:GgeE3}
\end{equation}

\textbf{Step 3:} We show that the kernel of the central charge of $\sigma$ is in Area I and is below the line through $E_1(-3)[2]$ and $E_3$, i.e. the open region bounded by line segments $\overline{RE_3}$, $\overline{E_3E_3^l}$ and $\overline{E_3^l R}$, with $R\coloneqq L_{E_1(-3)E_3}\cap L_{E_1E_2}$ as in Figure~\ref{fig:compareslope2}.

Let the central charge of $E_i[3-i]$ be $z_i$ for $i=1,2,3$. Let the object $G$ be of the form $E_1^{\oplus n_1}\rightarrow E_2^{\oplus n_2}\rightarrow E_3^{\oplus n_3}$.

By Lemma  \ref{lem:twotermstabobj} and a same argument as that in Lemma \ref{lem:lowbdforgldonlegs}, we know that the object $\mathsf{L}_{E_3}(E_1(-3)[2])$ is of the form $\Cone(E_1^{\oplus r_1}\rightarrow E_2^{\oplus r_2})$ in the heart $\tilde \cA$, and it is stable
with respect to every stability condition in $\Theta_\cE(\phi_2<\phi_1+1)$. By (\ref{eq:GgeE3}), we have $\Hom(G,E_3)=0$. By (\ref{eq:nonvanofe3fge1}), we have $\Hom(G,\mathsf{L}_{E_3}(E_1(-3)[2]))\neq 0$. Therefore, we have
\begin{equation}
    \phi(\mathsf{L}_{E_3}(E_1(-3)[2]))>\phi(G)>\phi(E_3).
    \label{eq:phaseofG}
\end{equation}
Therefore, the kernel of the central charge of $\sigma$ is in Area I and is below the line $L_{E_1(-3)E_3}$ as in Figure \ref{fig:compareslope2}.

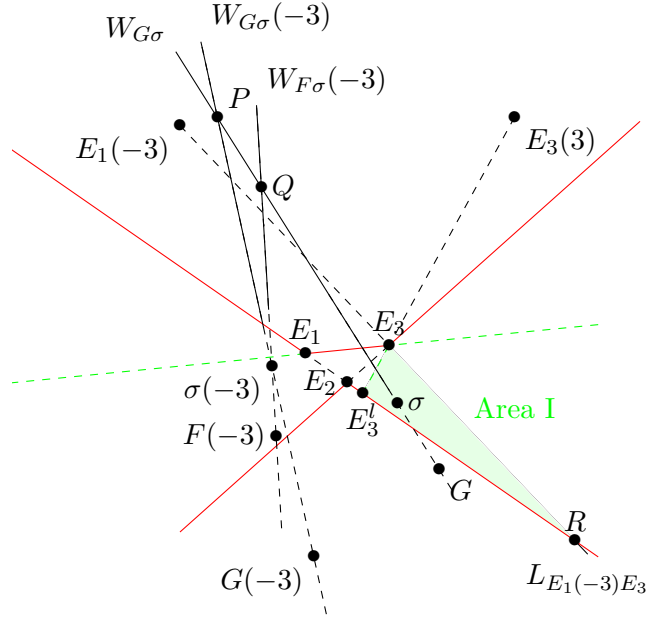
\begin{figure}[ht]\centering
\begin{tikzpicture}[yscale=.55,xscale=.55,>=stealth]
\clip (-8,-6.5) rectangle (8,9);
\tikzset{%
    add/.style args={#1 and #2}{
        to path={%
 ($(\tikztostart)!-#1!(\tikztotarget)$)--($(\tikztotarget)!-#2!(\tikztostart)$)%
  \tikztonodes},add/.default={.2 and .2}}
}

\newcommand\xg{2.2}
\newcommand\yg{-2.8}

\newcommand\xs{1.2}
\newcommand\ys{-1.2}

\coordinate (E1) at (-1,0);
\coordinate (E2) at (0,-0.7);
\coordinate (E3) at (1,0.2);
\coordinate (F1) at (-4,5.5); % F1=E_1(-3)
\coordinate (F3) at (4,5.7); % F3=E_3(3)
\coordinate (F) at (-1.7,-2); % F(-3)
\coordinate (G) at (\xg,\yg);
\coordinate (C) at (\xs,\ys); % sigma
\coordinate (SS) at (\xs-3,\ys-3*\xs+4.5); % sigma(-3)
\coordinate (G3) at (\xg-3,\yg-3*\xg+4.5); % G(-3)

\coordinate (Q) at (intersection of  {SS}--{F} and {G}--{C});
\coordinate (P) at (intersection of  {SS}--{G3} and {G}--{C});
\coordinate (R) at (intersection of  {F1}--{E3} and {E1}--{E2});

\draw [add= -1 and 0.21] (F3) to (E3) coordinate (EL3);
\draw [add= -1 and 0.95] (F1) to (E3) node [below] {$L_{E_1(-3)E_3}$};

\draw[white,fill=green!10] (R)--(E3)--(EL3)--cycle;

\draw[dashed] (F3) -- (E3);
\draw[dashed] (E2) -- (E3);
\draw[dashed] (E2) -- (E1);
\draw[dashed] (F1) -- (E3);

\draw[red] ($(E1)!-8!(E2)$)--(E1);
\draw[red] (E1) -- (E3);
\draw[red] (E3)--($(E3)!-6!(E2)$);
\draw[red] ($(E1)!7!(E2)$)--(E2);
\draw[red] (E2)--($(E3)!5!(E2)$);

\draw[dashed,green] (E1)--($(E3)!6!(E1)$);
\draw[dashed,green] ($(E1)!3.5!(E3)$)--(E3);
\draw[dashed,green] (E3)--(EL3);

%\draw [red] (-1.3,5) node {Area A};
%\draw [red] (0,-3.5) node {Area B};
\draw [green] (4,-1.3) node {Area I};
%\draw [green] (-6,1.3) node {Area II};

\draw [add=0.3 and 1.7, dashed] (G) to (C);
\draw [add=-1.2 and 5.3] (G) to (C)  node[above left] {$W_{G\sigma }$};
\draw [add=0.3 and 1.7, dashed] (G3) to (SS);
\draw [add=-1.2 and 1.7] (G3) to (SS)  node[above right] {$W_{G\sigma }(-3)$};
\draw [add=-2.3 and 1.3] (G3) to (SS);

\draw [add=1.3 and 3.7, dashed] (F) to (SS);
\draw [add=-1.8 and 3.7] (F) to (SS)  node[above right] {$W_{F\sigma}(-3)$};
\draw [add=-3 and 2] (F) to (SS);

\draw (G) node {$\bullet$} node [below right] {$G$};
\draw (C) node {$\bullet$} node [right] {$\sigma$};
\draw (G3) node {$\bullet$}node [below left] {$G(-3)$};
\draw (SS) node {$\bullet$}node [below left] {$\sigma(-3)$};

\draw (F) node {$\bullet$} node [left] {$F(-3)$};
\draw (E1) node {$\bullet$} node [above] {$E_1$};
\draw (E2) node {$\bullet$} node [left] {$E_2$};
\draw (E3) node {$\bullet$} node [above] {$E_3$};
\draw (F1) node {$\bullet$} node [below left] {$E_1(-3)$};
\draw (F3) node {$\bullet$} node [below right] {$E_3(3)$};
\draw (EL3) node {$\bullet$} node [below]{$E_3^l$};
\draw (R) node {$\bullet$} node [above] {$R$};
\draw (P) node {$\bullet$} node [above right] {$P$};
\draw (Q) node {$\bullet$} node [right] {$Q$};

cc\end{tikzpicture}
\caption{Comparing the phases of $G$ and $F(-3)$.}
\label{fig:compareslope2}
\end{figure}

\textbf{Step 4:} We show that the wall $W_{G\sigma}$ intersects the wall $W_{G\sigma}(-3)$. We denote the intersection point by $P\coloneqq W_{G\sigma}\cap W_{G\sigma}(-3)$.

Consider the line $L_{G\sigma}$ through the reduced character of $G$ and the kernel of the central charge of $\sigma$, which is in Area I. In particular, the stability condition  $$\sigma\in\Theta^+_\cE(\phi_2<\phi_1+1,\phi_3<\phi_1+2).$$ By (\ref{eq:phaseofG}), the line $L_{G\sigma}$ intersects the line segment $\overline{E_1(-3)E_3}$. Therefore, the line $L_{G\sigma}$ intersects the region $\operatorname{MZ}^c_\cE$.

By Lemma \ref{lem:nestedforalgstab}, the object $G$ is stable with respect to every stability condition in $\sigma$ in $\Theta^+_\cE(\phi_2<\phi_1+1,\phi_3<\phi_1+2)$ with kernel on $L_{G\sigma}$. Note that there exists a point $(1,s_0,q_0)$ in $\operatorname{MZ}^c_\cE\cap L_{G\sigma}$, by Lemma \ref{lemma:geoandalg}, the object $G$ is $\sigma_{s_0,q_0}$-stable.

Recall the wall $W_{G\sigma}\coloneqq\{(1,s,q)\in \cccp|$ the line segment along the line $L_{G\sigma}$ that is above the Le Potier curve $\mathrm{C}_{\mathrm{LP}} \}$.
By the Bertram's nest wall theorem \cite[Corollary 1.24]{LZ:mmpP2}, the object $G$ is $\sigma_{s,q}$-stable  for every $(1,s,q)$ on $W_{G\sigma}$. Note that $W_{G\sigma}$ intersects the line segment $\overline{E_1E_3}$, but does not intersects the line segment $\overline{E_1(-3)E_1}$ or $\overline{E_3E_3(3)}$, both of which are above the Le Potier curve $\mathrm{C}_{\mathrm{LP}}$. Therefore the horizental length of $W_{G\sigma}$ is greater than $3$ when $W_{G\sigma}$ is not the vertical wall. Let
$$W_{G\sigma}(-3)\coloneqq\{(1,s-3,q-3s+\frac{9}2)\; |\; (1,s,q)\in W_{G\sigma}\},$$
then $G(-3)$ is $\sigma_{s,q}$-stable  for every $(1,s,q)$ on $W_{G\sigma}(-3)$.
The wall $W_{G\sigma}(-3)$ intersects the wall $W_{G\sigma}$ at some point $P$ and
\begin{equation}
    \phi_{P}(G(-3))<\phi_{P}(G). \label{eq:comparGG3}
\end{equation}

As for the only exceptional case that $W_{G\sigma}$ is the vertical wall, we can view that the point $P$ is at $(0,0,1)$. This will not affect the statement in the next step.

\textbf{Step 5:} When $s(F)>s(E_3)$, we show that the wall $W_{F\sigma}(-3)$ intersects the wall $W_{G\sigma}$. We denote the intersection point by $Q\coloneqq W_{G\sigma}\cap W_{F\sigma}(-3)$.

By (\ref{eq:phaseofG}), we have the same bounds for $F$
\begin{equation}
    \phi(\mathsf{L}_{E_3}(E_1(-3)[2]))>\phi(G)>\phi(F)>\phi(E_3).
    \label{eq:phaseofF}
\end{equation}
The horizontal length of $W_{F\sigma}(-3)$ is greater than $3$ when it is not vertical. Note that the slope of $W_{F\sigma}(-3)$ is less than that of $W_{G\sigma}(-3)$, the segment $W_{F\sigma}(-3)$ intersects $W_{G\sigma}$ at $Q$ on the line segment $\overline{P\sigma}$. The fact that $\phi_{\sigma}(G)>\phi_{\sigma}(F)$ implies $\phi_{\sigma(-3)}(G(-3))>\phi_{\sigma(-3)}(F(-3))$. Both $F(-3)$ and $G$ are $\sigma_Q$ stable. We then compare their phases at $Q$ by using (\ref{eq:comparGG3}) as follows:
$$
\phi_Q(G)=\phi_P(G)> \phi_P(G(-3))=\phi_{\sigma(-3)}(G(-3))>\phi_{\sigma(-3)}(F(-3))=\phi_Q(F(-3)).
$$
So $\Hom(G,F(-3))=0$. By Serre duality, we have $\Hom(F,G[2])=0.$

\textbf{Step 6:} When $s(F) \leq s(E_3)$, we reduce this case to Proposition \ref{prop:mingldimpart} .

Note that $F$ is of the form $E_1^{\oplus a_1}\rightarrow E_2^{\oplus a_2}\rightarrow E_3^{\oplus a_3}$, we have $\Hom(E_3,F)\neq 0$ when $a_3\neq 0$. The object $F$ is either of the form $\Cone(E_1^{\oplus a_1}\rightarrow E_2^{\oplus a_2})[1]$ or $E_3$. Let $(1,s_0,q_0)$ be a point in $\operatorname{MZ}^c_\cE\cap L_{G\sigma}$. By Lemma \ref{lem:twotermstabobj}, in any case, $F$ is $\sigma_{s_0,q_0}$-stable. By Lemma \ref{lem:nestedforalgstab} and Lemma \ref{lemma:geoandalg}, the object $G$ is also  $\sigma_{s_0,q_0}$-stable and has phase
$$\phi_{s_0,q_0}(G)>\phi_{s_0,q_0}(F).$$
By Proposition \ref{prop:mingldimpart}, we have $\Hom(F,G[2])=0.$

As a summary, we have shown that $\Hom(F,G[2])=0$ when $\phi(F)<\phi(\mathsf{L}_{E_3}E_3(3))$ or $\phi(G)>\phi(E_3)$. In particular, we have $\gldim(\sigma)=\phi(E_3)-
    \phi(\mathsf{L}_{E_3}E_3(3))+2$.
\end{proof}

\begin{proof}[Proof for Proposition \ref{prop:glonthestab}]
When $\sigma\in \Theta_{\cE}\setminus \left(\Theta^{\mathrm{right}}_{E_1}\cup \Theta^{\mathrm{left}}_{E_3}\cup \Theta^{\mathrm{Pure}}_\cE\right)$, the global dimension is computed in Proposition \ref{prop:mingldimpart}.

When $\sigma \in \Theta^{\mathrm{left}}_{E_3}$, the global dimension is computed in Propositions \ref{prop:stabobjinfarleg} and  \ref{prop:glfunconhipleg}.

When $\sigma \in \Theta^{\mathrm{right}}_{E_1}$, we take the derived dual stability condition $\sigma^\vee\in \Theta^{\mathrm{left}}_{\cE^\vee,E_1^\vee}$, where $\cE^\vee$ is the dual exceptional triple $\{E_3^\vee,E_2^\vee,E_1^\vee\}$. We reduce to the previous case and have
\begin{align*}
\gldim(\sigma) & =\gldim(\sigma^\vee)=\phi^\vee(E^\vee_1)-\phi^\vee(\mathsf{L}_{E_1^\vee}(\Serre^{-1} E_1^\vee)) \\
&= -\phi(E_1)-\phi^\vee\left((\mathsf{R}_{E_1}(\Serre E_1))^\vee\right)=\phi(\mathsf{R}_{E_1}(\Serre E_1))-\phi_1.
\end{align*}

When $\sigma\in \Theta^{\mathrm{Pure}}_\cE$, by Lemma \ref{lemma:purepure}, the only stable objects are $E_i[m]$ for $E_i\in \cE$ and $m\in \mathbb Z$. As $\cE$ is a strong exceptional collection, we have $\Hom(E_i,E_j[m])\neq 0$ if and only $j\geq i$ and $m=0$. So the result is clear.
\end{proof}

\begin{remark}\label{rem:gldim3}
Following the notations in Remark \ref{rem:MZvsM}, for any exceptional bundle $E$, we associate two regions $\operatorname{MZ}^l_E$ and $\operatorname{MZ}^r_E$, which consist of geometric stability conditions. Moreover, if $\sigma\in \operatorname{MZ}^l_E \setminus \overline{E^lE}$ or $\sigma\in \operatorname{MZ}^r_E \setminus \overline{EE^r}$, we have $2<\gldim(\sigma)<3$.
\end{remark}

\begin{corollary}
The global dimension function 
$$
\gldim\colon\Stap(\PP^2)\to\RR_{\geq 0}
$$
has minimum value $2$ and $\gldim \Stap{\PP^2}=[2, \infty)$. Moreover, the subspace $\gldim^{-1}(2)$ is contained in $\overline{\Stab^{\Geo}(\PP^2)}$, and is contractible.
\label{cor:contra2}
\end{corollary}
\begin{proof}
The image of $\gldim$ follows from Proposition \ref{prop:gldiminParabola},  Proposition \ref{prop:glonthestab} and the description of $\Stap{\PP^2}$ (\ref{eq:stabmfd}). The contractibility of $\gldim^{-1}(2)$ is clear.
\end{proof}

%=========================================================
\section{Contractibility via global dimension}
%=========================================================
We denote by $\gldim^{-1}(I)$ by the space of all stability conditions in the component $\Stab^\dag(\PP^2)$ with global dimension in $I$ for an interval $I\subset\RR$. Based on Proposition \ref{prop:glonthestab} and the cell-decomposition description for $\Stab^\dag(\PP^2)$, our main result shows that the connected component $\Stab^\dag(\PP^2)$ is contractible via the global dimension:
\begin{theorem}\label{thm:3}
For any $x>2$, the space $\gldim^{-1}\left([2,x)\right)$ contracts to $\gldim^{-1}(2)$.
\end{theorem}
\begin{proof}
By Proposition  \ref{prop:glonthestab}, Remark \ref{rem:meandmze} and  \cite[Corollary 3.5, Theorem 3.9]{Li:spaceP2}, the space of preimage
$\gldim^{-1}\left([2,x)\right)$ has a cell decomposition as
\begin{equation*}
    \gldim^{-1}(2)\bigcup
\left(\bigsqcup_{E}\left( \Theta^{\text{left}}_{E}(x)\bigsqcup \Theta^{\text{right}}_{E}(x)\right) \bigsqcup \left(\bigsqcup_{\cE}
\Theta^{\mathrm{Pure}}_{\cE}(x)\right)\right),
\end{equation*}
where $E$ runs all exceptional bundles, and $\cE$ runs all exceptional triples, and the notation $\Theta_*^\dag(x)$ stands for $\Theta_*^\dag\cap \gldim^{-1}\left([2,x)\right)$. 

By Proposition \ref{prop:glonthestab}, we have $\Theta^{\mathrm{Pure}}_{\cE}(x)=\Theta^{\mathrm{Pure}}_{\cE}(\phi_3-\phi_1<x)$. Each $\Theta^{\mathrm{Pure}}_{\cE}(x)$ has an open neighborhood, say, $\Theta_{\cE}(\phi_3-\phi_2>\frac{1}{2}, \phi_2-\phi_1>\frac{1}{2}, \phi_3-\phi_1<x)$, in
$\Theta_{\cE}(x)$ which does not intersect any other
$\Theta^{\mathrm{Pure}}_{\cE'}(x)$.  As $\Stab^\dag(\PP^2)$ admits a
metric, we may then choose open neighborhoods of
$\Theta^{\mathrm{Pure}}_{\cE}(x)$'s which do not intersect with each other.
By the cell decomposition, the space $\gldim^{-1}\left([2,x)\right)$ contracts to its subspace $$A(x)\coloneqq\gldim^{-1}(2)\bigcup
\left(\bigsqcup_{E\text{ exceptional bundles}}\left( \Theta^{\text{left}}_{E}(x)\bigsqcup \Theta^{\text{right}}_{E}(x)\right)\right).$$

For each exceptional object $E$, let $\cE=\{E_1,E_2,E_3\}$ be an exceptional collection such that $E_3=E$. By Proposition \ref{prop:glonthestab} and Lemma \ref{lem:basicexc}, we have
\begin{align*}
\Theta^{\text{left}}_{E}(x)\cong & \big\{(m_1,m_2,m_3,\phi_1,\phi_2,\phi_3)\in (\mathbb R_{>0})^3\times
\mathbb R^3|\phi_1<\phi_2<\phi_1+1, \\ & \phi_3>\phi_2+1+\frac{1}\pi \arctan\left(\frac{\sin((\phi_1+1-\phi_2)\pi)}{\cos((\phi_1+1-\phi_2)\pi)+\frac{m_2}{m_1}h}\right)>\phi_3-x+2\big\},
\end{align*}
where $h=\mathrm{hom}(E_1,E_2)-\frac{\mathrm{hom}(E_2,E_3)}{\mathrm{hom}(E_1,E_3)}$. Therefore, the space $\Theta^{\text{left}}_{E}(x)$ contracts to  $\Theta^{\text{left}}_{E}(x)\cap \gldim^{-1}(2)$.

By Remark \ref{rem:meandmze} and \cite[Lemma 3.7 ]{Li:spaceP2}, each $\Theta^{\text{left}}_{E}(x)$ has an open
neighborhood in $A(x)$, which does not intersect any other
$\Theta^{\text{left}}_{E'}(x)$ or $\Theta^{\text{right}}_{E'}(x)$. Same argument works for all $\Theta^{\text{right}}_{E}(x)$, we may therefore contract all $\Theta^{\text{left}}_{E}(x)$ and $\Theta^{\text{right}}_{E}(x)$ in $A(x)$ simultaneously to $\gldim^{-1}(2)$, which is a contractible space.
\end{proof}

%=========================================================
\section{Inducing stability conditions from projective plane to the local projective plane}\label{sec:inducing}
%=========================================================
Let $Y$ be the total space of the canonical bundle of $\PP^2$,
and $i:\PP^2\hookrightarrow Y$ be the inclusion of the zero-section.
We write $\D^b_{\PP^2}(Y)$ for the subcategory of $\D^b(Y)$ of complexes with bounded cohomology, such that all of its cohomology sheaves are supported on the zero-section.
The space of Bridgeland stability conditions on $\D^b_{\PP^2}(Y)$ has been studied by Bayer and Macr{\`{\i}} \cite{BM11}.
In this section, we prove that the stability conditions in
$\gldim^{-1}(2)\subset\overline{\Stab^{\Geo}(\PP^2)}$
can be used to induce stability conditions on $\D^b_{\PP^2}(Y)$
by Ikeda-Qiu's inducing theorem, via $q$-stability conditions on Calabi--Yau-$\XX$ categories.

Following the notion in \cite{IQ18a},
we have the Calabi--Yau-$\XX$ version of $\D_\infty(\PP^2)$
\begin{gather}\label{eq:DX P2}
    \D_\XX(\PP^2)\coloneqq\D^b_{c,\CC^*}(Y).
\end{gather}
By \cite[Proposition~3.14]{IQ18a}, we have
$\D_\XX(\PP^2)\cong\D_{\mathrm{fd}}(\Gamma_\XX(\widetilde{Q}_{\mathrm{gr}},W_{\mathrm{gr}}))$
with $\ZZ\oplus\ZZ[\XX]$ graded quiver $\widetilde{Q}_{\mathrm{gr}}$ as follows and potential $W_{\mathrm{gr}}=\sum_{i=1}^3 (x_iy_iz_i-x_iz_iy_i)$,
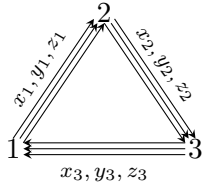
\begin{center}
\begin{tikzpicture}[scale=0.6,
  arrow/.style={->,>=stealth},
  equalto/.style={double,double distance=2pt},
  mapto/.style={|->}]
\node (x2) at (0,2.){};
\node (x1) at (2,-1){};
\node (x3) at (-2,-1){};
  \node at (x1){$3$};
  \node at (x2){$2$};
  \node at (x3){$1$};
\foreach \n/\m in {1/3}
    {\draw[arrow] (x\n.150) to (x\m.30);
     \draw[arrow] (x\n.-150) to (x\m.-30);
     \draw[arrow] (x\n) to (x\m);}
\foreach \n/\m in {2/1}
    {\draw[arrow] (x\n.150+120) to (x\m.30+120);
     \draw[arrow] (x\n.-150+120) to (x\m.-30+120);
     \draw[arrow] (x\n) to (x\m);}
\foreach \n/\m in {3/2}
    {\draw[arrow] (x\n.150-120) to (x\m.30-120);
     \draw[arrow] (x\n.-150-120) to (x\m.-30-120);
     \draw[arrow] (x\n) to (x\m);}
\draw[font=\scriptsize](-.8,1.8)node[left,rotate=60]{$x_1,y_1,z_1$}
    (.6,1.8)node[right,rotate=-54]{$x_2,y_2,z_2$}
    (0,-1.2)node[below]{$x_3,y_3,z_3$};
\end{tikzpicture}
\end{center}
where $\deg x_3, y_3, z_3=3-\XX$ and gradings of other arrows are zero.
Here $\Gamma_\XX(\widetilde{Q}_{\mathrm{gr}},W_{\mathrm{gr}})$ is the Calabi--Yau-$\XX$ Ginzburg dg algebra \cite{IQ18a,IQ18b}.
Note that there is a canonical fully faithful embedding
\[
    \D_\infty(\PP^2)\to\D_\XX(\PP^2)
\]
whose image is an $\XX$-baric heart of $\D_\XX(\PP^2)$ in the sense of \cite[Definition~2.17]{IQ18a}.

Finally, we have the 3-reduction of $\D_\XX(\PP^2)$ (see \cite[Example~3.16]{IQ18a})
\begin{gather}\label{eq:D3 P2}
    \D_3(\PP^2)\coloneqq\D_\XX(\PP^2)\sslash[\XX-3]\cong\D^b_{\PP^2}(Y)
\end{gather}
which is equivalent to the derived category of coherent sheaves on the local $\PP^2$.

\begin{lemma} \label{Lem:induce}
Consider the composition of functors in the inducing process
\begin{multline*}
\Phi:\D_\infty(\PP^2) \xrightarrow{\sim} \D^b(\k Q/R) \to
\D_{\mathrm{fd}}(\Gamma_\XX(\widetilde{Q}_{\mathrm{gr}},W_{\mathrm{gr}})) \to  \\
\D_{\mathrm{fd}}(\Gamma_\XX(\widetilde{Q}_{\mathrm{gr}},W_{\mathrm{gr}})) \sslash [\XX-3] \xrightarrow{\sim}
\D^b(\mod-J(\widetilde{Q},W)) \xrightarrow{\sim} \D^b_{\PP^2}(Y).
\end{multline*}
Then
$\Phi = i_*: \D_\infty(\PP^2) \to\D^b_{\PP^2}(Y).
$
\end{lemma}

\begin{proof}
Let $E_i=\O_{\PP^2}(i)$.
Then the first  equivalence $\D_\infty(\PP^2) \cong \D^b(\k Q/R) $ in $\Phi$ is given by
\[
\Hom^\bullet(\bigoplus_{i=0}^2 E_i, -): \D_\infty(\PP^2) \to  \D^b(\k Q/R),
\]
and the last equivalence $\D^b_{\PP^2}(Y) \cong\D^b(\mod-J(\widetilde{Q},W))$ in $\Phi$ is given by
\[
\Hom^\bullet(\bigoplus_{i=0}^2 \pi^*E_i, -): \D^b_{\PP^2}(Y) \to  \D^b(\mod-J(\widetilde{Q},W)),
\]
where $\pi:Y\to\PP^2$ is the projection \cite{Bri05}.
The lemma then follows from
$$
\Hom^\bullet(\pi^*\mathcal{E}, i_*{\mathcal{F}}) =
\Hom^\bullet(\mathcal{E}, \pi_*i_*\mathcal{F}) =
\Hom^\bullet(\mathcal{E}, \mathcal{F}).
$$
\end{proof}

Now we recall the inducing construction of stability conditions
from the projective plane to the local projective plane,
through the `$q$-stability conditions' introduced by Ikeda and Qiu \cite{IQ18a}.

\begin{construction}
Let $\sigma_\infty=(Z_\infty,\sli_\infty)$ be a stability condition in
$\gldim^{-1}(2)\subset\overline{\Stab^{\Geo}(\PP^2)}$.
\begin{itemize}
  \item By \cite[Theorem.~2.25]{IQ18a},
  there is an induced $q$-stability conditions $(\sigma,s)$ in $\QStab\D_\XX(\PP^2)$
  with parameter $s=3$,
  as constructed in \cite[Cons.~2.18]{IQ18a}.
  \item By \cite[Theorem.~2.16]{IQ18a}, $(\sigma,s)$ projects to a stability condition $\sigma_3$
  in the principal (connected) component $\Stap\D_3(\PP^2)$.
\end{itemize}
Denote by
\[
    \iota_3\colon \gldim^{-1}(2)\to \Stap \D_3(\PP^2)
\]
the map of the above inducing process.
\end{construction}

\begin{proposition} \label{Prop:induce}
The inducing map $\iota_3$ is injective.
Moreover, it factors through the isomorphism between the spaces of geometric stability conditions on $\PP^2$ and local $\PP^2$:
\[
    \iota_3\colon \gldim^{-1}(2)\hookrightarrow\overline{\Stab^{\Geo}(\PP^2)} \iso \overline{\Stab^{\Geo} (\D^b_{\PP^2}  (Y)   ) } \hookrightarrow \Stap \D_3(\PP^2).
\]
\end{proposition}

\begin{proof}
A stability condition $\sigma$ on $\D^b_{\PP^2}(Y)$ is called \emph{geometric} if all skyscraper sheaves $i_*\cO_x$ of closed points $x\in\PP^2$ are $\sigma$-stable of the same phase.
By Lemma \ref{Lem:induce}, the inducing map $\iota_3$ maps geometric stability conditions on $\PP^2$ with global dimension 2 to geometric stability conditions on local $\PP^2$.

Let $\sigma=(Z,P)\in\gldim^{-1}(2)$ and
$\iota_3(\sigma)=(\widetilde{Z}, \widetilde{P})\in\Stab^{\Geo} (\D^b_{\PP^2}  (Y)   )$.
By Lemma \ref{Lem:induce}, we have
$Z = \widetilde{Z} \circ [i_*]
$,
where $[i_*]: K_0(\D_\infty(\PP^2)) \iso K_0(\D^b_{\PP^2}(Y))$.
By \cite[Theorem 2.5]{BM11} and \cite[Proposition 1.12]{Li:spaceP2},
any geometric stability condition on $\PP^2$ or local $\PP^2$ is uniquely determined by its central charge.
Moreover, the open set  $U\subset\Hom(K_0(\D_\infty(\PP^2)), \mathbb{C})$ consists of central charges of geometric stability conditions on $\PP^2$ and the open set
$\widetilde{U}\subset \Hom(K_0(\D^b_{\PP^2}(Y)), \mathbb{C})$ of central charges of geometric stability conditions on local $\PP^2$ coincide via the isomorphism $[i_*]$.
This proves the proposition.
\end{proof}

Finally, we remark that the whole connected component of stability conditions in $\Stap \D_3(\PP^2)$ can be obtained by inducing from stability conditions on $\D_\infty(\PP^2)$ and autoequivalences, since the translates of $\overline{\Stab^{\Geo} (\D^b_{\PP^2}(Y))}$
under the group of autoequivalences cover the whole connected component $\Stap \D_3(\PP^2) $ \cite[Theorem 1]{BM11}.

%If there are any kind of inducing, in the sense of stable objects in $\D_\infty$ remains stable in $D_3$, then $\gldim\leq 3$ due to Calabi--Yau-3 property. }

%=========================================================

\bibliography{refsAMS}{}
\bibliographystyle{alpha}

%=========================================================
\end{document}